\documentclass[11pt,a4paper,leqno]{amsart}
\usepackage{amssymb}
\newif\ifaddpics\addpicstrue
\ifaddpics\usepackage{graphicx}\usepackage{diagrams}\fi
\swapnumbers
\theoremstyle{plain}

\newtheorem{thm}{Theorem}[section]
\newtheorem{prop}[thm]{Proposition}
\newtheorem{cor}[thm]{Corollary}
\newtheorem{lemma}[thm]{Lemma}

\theoremstyle{definition}
\newtheorem{defn}[thm]{Definition}

\theoremstyle{remark}

\newcommand{\proofof}[1]{\end{#1}\begin{proof}} 
\newcommand{\emphdef}{\textit}                  
\newcommand{\tcite}[1]{\textup{\cite{#1}}}      
\numberwithin{equation}{section}

\renewcommand{\labelenumi}{\textup{(\roman{enumi})}}

\newcommand{\acknowledge}{\subsection*{Acknowledgements}}
\newcommand{\thismonth}{\ifcase\month\or
  January\or February\or March\or April\or May\or June\or
  July\or August\or September\or October\or November\or December\fi
  \space\number\year}
\makeatletter
\newcommand{\low}{\@ifnextchar^{}{^{\vphantom x}}}
\newcommand{\high}{\@ifnextchar_{}{_{\vphantom I}}}
\makeatother
\DeclareSymbolFont{script}{U}{eus}{m}{n}
\DeclareSymbolFontAlphabet{\mathscr}{script}
\DeclareMathSymbol{\EuWedge}{0}{script}{"5E}
\DeclareMathAlphabet{\mathrmsl}{OT1}{cmr}{m}{sl}
\newcommand{\rssymb}[2]{\newcommand{#1}{{\mathrmsl{#2}}}}
\newcommand{\calsymb}[2]{\newcommand{#1}{{\mathcal{#2}}}}
\newcommand{\bbsymb}[2]{\newcommand{#1}{{\mathbb{#2}}}}
\newcommand{\lieoper}[2]{\newcommand{#1}{\mathop
  {\mathfrak{#2}\null}\nolimits}}
\newcommand{\oper}[3][n]{\newcommand{#2}{\mathop
  {\mathrm{#3}\null}\ifx n#1\nolimits\else\limits\fi}}
\newcommand{\rsoper}[3][n]{\newcommand{#2}{\mathop
  {\mathrmsl{#3}\null}\ifx n#1\nolimits\else\limits\fi}}
\bbsymb\C{C} \bbsymb\F{F} \bbsymb\HQ{H}\bbsymb\N{N} \bbsymb\Q{Q}
\bbsymb\R{R} \bbsymb\U{U} \bbsymb\V{V} \bbsymb\W{W} \bbsymb\Z{Z}
\calsymb\cA{A} \calsymb\cB{B} \calsymb\cC{C} \calsymb\cD{D} \calsymb\cE{E}
\calsymb\cF{F} \calsymb\cG{G} \calsymb\cH{H} \calsymb\cI{I} \calsymb\cJ{J}
\calsymb\cK{K} \calsymb\cL{L} \calsymb\cM{M} \calsymb\cN{N} \calsymb\cO{O}
\calsymb\cP{P} \calsymb\cQ{Q} \calsymb\cR{R} \calsymb\cS{S} \calsymb\cT{T}
\calsymb\cU{U} \calsymb\cV{V} \calsymb\cW{W} \calsymb\cX{X} \calsymb\cY{Y}
\calsymb\cZ{Z}
\newcommand{\eps}{\varepsilon}
\newcommand{\gam}{\gamma} 

\renewcommand{\geq}{\geqslant} \renewcommand{\leq}{\leqslant}
\oper\End{End} \oper\Hom{Hom}                    
\oper\Sym{Sym} \oper\Skew{Skew}
\oper\Aut{Aut}                                   
\oper\GL{GL} \oper\SL{SL}\oper\Symp{Sp}
\oper\CO{CO} \oper\On{O} \oper\SO{SO} \oper\Pin{Pin} \oper\Spin{Spin}
\oper\CU{CU} \oper\Un{U} \oper\SU{SU}
\rsoper\Diff{Diff} \rsoper\SDiff{SDiff}
\lieoper\der{der}                                
\lieoper\gl{gl} \lieoper\sgl{sl}\lieoper\symp{sp}
\lieoper\co{co} \lieoper\so{so} \lieoper\spin{spin}
\lieoper\cu{cu} \lieoper\un{u}  \lieoper\su{su}
\rsoper\Vect{Vect} \rsoper\Ham{Ham}
\newcommand{\ip}[1]{\langle#1\rangle}

\newcommand{\norm}[2][]{|\mkern-2mu|#2|\mkern-2mu|
  _{\lower1pt\hbox{${}_{#1}$}}}
\newcommand{\Norm}[2][]{\bigl|\mkern-3mu\bigr|#2\bigr|\mkern-3mu\bigr|
  _{\lower1pt\hbox{${}_{#1}$}}}

\newcommand{\lie}[1]{\mathfrak{#1}}

\newcommand{\punc}[1]{\smallsetminus\{#1\}}
\newcommand{\restr}[1]{|_{#1}\low}

\newcommand{\cross}{\mathbin{{\times}\!}\low}
\newcommand{\dual}{^{*\!}}
\newcommand{\mult}{^{\scriptscriptstyle\times}}

\newcommand{\act}{\mathinner\cdot}          
\newcommand{\dsum}{\oplus}                  
\newcommand{\tens}{\otimes}                 
\newcommand{\Wedge}{\EuWedge}               
\newcommand{\skwend}{\mathinner\vartriangle}
\newcommand{\setdif}{\smallsetminus}
\newcommand{\connect}{\#}                   
\newcommand{\del}{\partial}                 
\newcommand{\dbar}{\overline\partial}       
\newcommand{\Proj}{\mathrmsl{P}}            
\newcommand{\RP}[1]{\R\Proj^{#1}}           
\newcommand{\CP}[1]{\C\Proj^{#1}}           
\newcommand{\HP}[1]{\HQ\Proj^{#1}}          
\newcommand{\Cinf}{\mathrm{C}^\infty}       
\rsoper\dimn{dim}                           
\rsoper\grad{grad}                          
\rsoper\kernel{ker}\rsoper\image{im}        
\rsoper\alt{alt}   \rsoper\sym{sym}         
\rsoper\Ad{Ad}     \rsoper\ad{ad}           
\rsoper\CoAd{CoAd} \rsoper\coad{coad}       
\rsoper\trace{tr}  \rsoper\trfree{tf}       
\rsoper\detm{det}                           
\rsoper\Vol{Vol}                            
\rsoper\divg{div}                           
\rsoper\sign{sign}
\rssymb\iden{id}                            
\rssymb\vol{vol}                            
\oper\Imag{Im}
\newcommand{\sd}{{\raise1pt\hbox{$\scriptscriptstyle +$}}}
\newcommand{\asd}{{\raise1pt\hbox{$\scriptscriptstyle -$}}}
\newcommand{\sdasd}{{\raise1pt\hbox{$\scriptscriptstyle\pm$}}}
\newcommand{\asdsd}{{\raise1pt\hbox{$\scriptscriptstyle\mp$}}}
\newcommand{\gmw}{\mathrmsl w}

\newcommand{\iI}{{\boldsymbol i}}
\newcommand{\jJ}{{\boldsymbol j}}
\newcommand{\kK}{{\boldsymbol k}}
\newcommand{\ratio}{\mathbin:}

\rsoper\scal{scal}
\def\kahl/{k{\"a}hler}
\begin{document}
\title{Selfdual Einstein metrics with torus symmetry}
\author{David M. J. Calderbank}
\address{Department of Mathematics and Statistics\\
University of Edinburgh\\ King's Buildings, Mayfield Road\\
Edinburgh EH9 3JZ\\ Scotland.}
\email{davidmjc@maths.ed.ac.uk}
\author{Henrik Pedersen}
\address{Department of Mathematics and Computer Science\\
SDU-Odense University\\ Campusvej 55\\ DK-5230 Odense M\\ Denmark.}
\email{henrik@imada.sdu.dk}
\date{May 2000}
\begin{abstract}
It is well known that any $4$-dimensional hyper\kahl/ metric with two
commuting Killing fields may be obtained explicitly, via the Gibbons--Hawking
Ansatz, from a harmonic function invariant under a Killing field on $\R^3$.
In this paper, we find all selfdual Einstein metrics of nonzero scalar
curvature with two commuting Killing fields. They are given explicitly in
terms of a local eigenfunction of the Laplacian on the hyperbolic plane.  We
discuss the relation of this construction to a class of selfdual spaces found
by Joyce, and some Einstein--Weyl spaces found by Ward, and then show that
certain `multipole' hyperbolic eigenfunctions yield explicit formulae for the
quaternion-\kahl/ quotients of $\HP{m-1}$ by an $(m-2)$-torus first studied by
Galicki and Lawson.  As a consequence we are able to place the well-known
cohomogeneity one metrics, the quaternion-\kahl/ quotients of $\HP{2}$ (and
noncompact analogues), and the more recently studied selfdual Einstein
Hermitian metrics in a unified framework, and give new complete examples.
\end{abstract}
\maketitle

\section{Introduction}

We present in this paper an explicit classification of $4$-dimensional
Einstein metrics with selfdual Weyl curvature and two linearly independent
commuting Killing fields. We refer to these metrics as selfdual Einstein
metrics with torus symmetry, since they are the local form (on a dense open
set) of such metrics with an action of $T^2$, $S^1\cross\R$ or $\R^2$ by
isometries.

When the selfdual Einstein metric $g$ is scalar-flat, it is well known that
$g$ is locally hyper\kahl/ and that some linear combination of the two Killing
fields is triholomorphic---hence $g$ is determined by a harmonic function on
$\R^3$, via the Gibbons--Hawking Ansatz~\cite{GiHa:gmi}, and this harmonic
function is invariant under the transrotation of $\R^3$ coming from the second
Killing field of $g$. Therefore we focus on the case that the selfdual
Einstein metric has nonzero scalar curvature.

\begin{thm}\label{main} Let $F(\rho,\eta)$ be a solution
of the linear differential equation
\begin{equation*}
F_{\rho\rho}+F_{\eta\eta}=\frac{3F}{4\rho^2}
\end{equation*}
on some open subset of the half-space $\rho>0$, and consider the metric
$g(\rho,\eta,\phi,\psi)$ given by
\begin{equation}\label{metric}\begin{split}
g&=\frac{ F^2 - 4\rho^2(F_\rho^2 + F_\eta^2) }{4 F^2}\;
\frac{d\rho^2 + d\eta^2}{\rho^2}\\
&\qquad\qquad+\frac{
\bigl( (F - 2 \rho F_\rho) \alpha - 2 \rho F_\eta \beta \bigr)^2 +
\bigl( -2\rho F_\eta  \alpha + (F + 2 \rho F_\rho )\beta \bigr)^2
}
{ F^2\bigl(F^2-4\rho^2(F_\rho^2 + F_\eta^2)\bigr)},
\end{split}\end{equation}
where $\alpha=\sqrt\rho\,d\phi$ and $\beta=(d\psi+\eta\, d\phi)/\sqrt\rho$.
Then:
\begin{enumerate}
\item On the open set where $F^2>4\rho^2(F_\rho^2 + F_\eta^2)$, $g$ is a
selfdual Einstein metric of positive scalar curvature, whereas on the open set
where $0<F^2<4\rho^2(F_\rho^2 + F_\eta^2)$, $-g$ is a selfdual Einstein metric
of negative scalar curvature.

\item Any selfdual Einstein metric of nonzero scalar curvature with two
linearly independent commuting Killing fields is arises locally in this way
\textup(i.e., in a neighbourhood of any point, it is of the
form~\eqref{metric} up to a constant multiple\textup).
\end{enumerate}
\end{thm}

The metric $g$ is sufficiently explicit to make it straightforward, though
tedious, to check that it is selfdual and Einstein, for instance, by computing
the curvature of the Levi-Civita connection of $g$ on antiselfdual $2$-forms
using the Cartan calculus~\cite{Tod:com}. Hence the heart of the above Theorem
is part (ii), so we concentrate on this and explain why all selfdual Einstein
metrics with torus symmetry are of this form. The proof of (ii) will in fact
encode the selfdual Einstein condition (using work of Tod~\cite{Tod:sut}),
thus proving (i) at the same time.

There are three features of these metrics which need to be explained.  First,
and most remarkable, is the fact that the equation for $F$ is linear---this
means that we can ``superpose'' any two such metrics to yield a third.

Second, the equation for $F$ means that it is a local eigenfunction, with
eigenvalue $3/4$, of the Laplacian $\rho^{-2}(\del_\rho^2+\del_\eta^2)$ of
the hyperbolic metric $(d\rho^2+d\eta^2)/\rho^2$---in other words $F$
satisfies a natural differential equation on the hyperbolic plane.

Third, the level surfaces of constant $\phi,\psi$ are orthogonal to both
Killing fields $\partial/\partial\phi$ and $\partial/\partial\psi$---hence the
orbits of the induced $2$-dimensional local symmetry group are
\emph{surface-orthogonal}, i.e., the orthogonal distribution to the orbits is
integrable.

These three features are closely related. We originally found an explicit form
for selfdual Einstein metrics with torus symmetry by noticing that the local
quotient of such a metric by one of its Killing fields must be an
\emphdef{Einstein--Weyl space with an axial symmetry}, i.e., one of the spaces
found by Ward~\cite{Ward:sut} and studied in~\cite{Cal:gte}.  These
Einstein--Weyl spaces are given explicitly in terms of an
\emphdef{axisymmetric harmonic function} (AHF) on $\R^3$, i.e., a solution of
a linear differential equation. The geometry of the hyperbolic plane $\cH^2$
enters the picture because $\R^3\setdif\R$ is conformal to $\cH^2\cross
S^1$. Finally, these Einstein--Weyl spaces are the quotients of the conformal
metrics found by Joyce~\cite{Joy:esd}, who obtained a classification of
selfdual $4$-manifolds with two commuting \emph{surface-orthogonal} conformal
vector fields.

This original argument had the advantage of being an exercise in pure thought,
i.e., a combination of known results with no new computations.  However, it
led to a rather awkward description of the conformally Einstein metrics among
Joyce's selfdual spaces.

Therefore, we shall present a different, more self-contained proof of
Theorem~\ref{main}. Indeed, we proceed rather in reverse, by establishing
directly, in section~\ref{so}, the surface-orthogonality of the Killing
fields. To do this, we first review, in section~\ref{sdekv}, an isomorphism
between Killing fields and \emphdef{twistors} which plays an important role
throughout the paper---these twistors are essentially the same thing as
compatible scalar-flat \kahl/ metrics~\cite{Pon:ths}. We discuss this in
section~\ref{sfktoda}, together with Tod's description~\cite{Tod:sut} of
selfdual Einstein metrics with a Killing field in terms of the $\SU(\infty)$
Toda equation (see also~\cite{Prz:kvf,FiSa:hek}).

Returning to the general argument, the surface-orthogonality of the Killing
fields shows that the conformal structure of a selfdual Einstein metric with
torus symmetry is a `Joyce space'~\cite{Joy:esd}. We review the theory of
Joyce spaces and their quotient Einstein--Weyl spaces in section~\ref{jaew}
and explain how they are determined by two solutions of a spinor equation on
$\cH^2$ which is equivalent to the equation for AHFs. On the other hand, Tod's
analysis of selfdual Einstein metrics with a Killing field shows that these
two AHFs must be related in a special way. The essential idea is that the AHFs
are constructed from the eigenfunction $F$ by a B\"acklund transformation, and
a $2$-dimensional family is obtained because this transformation is coordinate
dependent. In sections~\ref{sghp} and~\ref{ttp} we use the spin geometry of
the hyperbolic plane to give a more natural description, and this allows us to
complete the proof of Theorem~\ref{main}.

In section~\ref{swann} we discuss the Swann bundle~\cite{Swa:hkqk} of our
selfdual Einstein metrics. This hyper\kahl/ $8$-manifold with two commuting
triholomorphic Killing fields is given locally by the generalized
Gibbons--Hawking construction of~\cite{HKLR:hkms,PePo:hkm} and we establish
the relationship between this construction and the eigenfunction $F$. Then,
comparing the Swann bundle with the hyper\kahl/ quotients of $\HQ^m$ by an
$(m-2)$-torus (see~\cite{BiDa:gtt}), we are led, in section~\ref{qkhk}, to
define `multipole' hyperbolic eigenfunctions and prove the following theorem.

\begin{thm} The selfdual Einstein metrics arising as
quaternion-\kahl/ quotients of quaternionic projective space $\HP{m-1}$ by an
$(m-2)$-dimensional family of commuting Killing fields are exactly the
metrics given by~\eqref{metric} where $F$ is a `positive $m$-pole solution' of
the hyperbolic eigenfunction equation on $\cH^2$.
\end{thm}
A more precise statement is given in Theorem~\ref{qkposthm}. In particular,
this characterizes the $m$-pole solutions corresponding to quaternion-\kahl/
quotients of $\HP{m-1}$ by an $(m-2)$-torus, which yield the compact selfdual
Einstein orbifold metrics of Galicki--Lawson~\cite{GaLa:qro} when $m=3$, and
Boyer--Galicki--Mann--Rees~\cite{BGMR:3s7} in general.

Replacing $\HQ^m$ with $\HQ^{p,q}$ leads to non-compact analogues of these
metrics, among which there are examples of complete selfdual Einstein metrics
of negative scalar curvature. Some of these metrics are well-known---and our
approach provides a unified description of them---but we also obtain new
examples.

\acknowledge The first author would like to thank Vestislav Apostolov, Roger
Bielawski, Paul Gauduchon, Kris Galicki, Elmer Rees, and Michael Singer for
stimulating discussions, and also the EPSRC, the Leverhulme Trust and the
William Gordon Seggie Brown Trust for the Research Fellowships that enabled
him to carry out this research.  Both authors are members of EDGE, Research
Training Network HPRN-CT-2000-00101, supported by the European Human Potential
Programme

\section{Selfdual Einstein metrics, twistors and Killing fields}\label{sdekv}

We begin by reviewing the relation between Killing fields and antiselfdual
twistors on a selfdual Einstein manifold, following~\cite{Bes:em,Tod:sut}.
Recall that a vector field $K$ is a Killing field of a metric $g$ if and only
if its covariant derivative $D^gK$ is a skew endomorphism of the tangent
bundle (the endomorphism is defined by $X\mapsto D^g_XK$).

In four dimensions skew endomorphisms decompose into selfdual and antiselfdual
parts: $\so(TM)=\so\low_\sd(TM)\dsum\so\low_\asd(TM)$. This is related (via
the metric) to the decomposition $\Wedge^2T\dual M=\Wedge^2_\sd T\dual
M\dsum\Wedge^2_\asd T\dual M$ of the bundle of $2$-forms into eigenspaces of
the Hodge $*$ operator. We denote the decomposition of sections of these
bundles by $A=A^\sd+A^\asd$.

An antiselfdual endomorphism $\Psi$ is called a \emphdef{twistor} if there is
a $1$-form $\gam$ such that $D^g_X\Psi=(\gam\skwend X)^\asd$ for all vector
fields $X$---here, for vector fields $X,Y$, we define $\gam\skwend X
(Y)=\gam(Y)X-\ip{X,Y}\sharp\gam$. It follows by taking a trace that $\gam$ is
a multiple of the divergence $\delta^g\Psi$ and so this is a linear
differential equation on $\Psi$, called the \emphdef{twistor equation}.

For a selfdual Einstein metric, the curvature $R^g\in\Cinf(M,\Wedge^2 T\dual
M\tens \so(TM))$ has only two irreducible components: the selfdual Weyl
curvature $W^\sd\in\Cinf(M,\Wedge^2_\sd T\dual M\tens \so\low_\sd(TM))$, and
the (normalized) scalar curvature $s^g = \tfrac16 \scal^g$.

This has strong consequences for Killing fields and twistors.

\begin{prop}\label{sdeK} Let $g$ be a selfdual Einstein metric.
\begin{enumerate}
\item Suppose $K$ is a Killing field of $g$ and let $\Psi=(D^gK)^\asd$.  Then
$D^g_X\Psi = \frac12s^g\bigl(\ip{K,\cdot}\skwend X\bigr)^\asd$ and so $\Psi$
is a twistor.
\item Suppose $\Psi$ is a twistor with $D^g_X\Psi=(\gam\skwend X)^\asd$.  Then
the dual vector field $\sharp^g\gam$ is a Killing field of $g$ and
$(D^g\gam)^\asd= \frac12s^g \ip{\Psi\cdot,\cdot}$.
\end{enumerate}
\proofof{prop} (i) Since $K$ is Killing and $g$ is selfdual Einstein, we have
\begin{equation*}
D^g_X D^g K = R^g_{X,K} = W^\sd_{X,K} + \tfrac12 s^g\ip{K,\cdot}\skwend X.
\end{equation*}
The antiselfdual part of this is what we want.

\noindent(ii) Differentiating the twistor equation again and skew
symmetrizing gives
\begin{equation*}
(D^g_X\gam\skwend Y)^\asd - (D^g_Y\gam\skwend X)^\asd =
[R^g_{X,Y},\Psi]=-\frac12 s^g[\ip{X,\cdot}\skwend Y,\Psi].
\end{equation*}
Contracting with another vector field $Z$ and taking the trace
over $X$ and $Z$, we obtain
\begin{equation*}
\sym D^g\gam+\tfrac12\delta^g\gam+(d\gam)^\asd
=s^g\ip{\Psi\cdot,\cdot}.
\end{equation*}
The right hand side is skew, so $D^g\gam$ is skew, $\sharp^g\gam$ is a Killing
field and $2(D^g\gam)^\asd =s^g\ip{\Psi\cdot,\cdot}$.
\end{proof}

\begin{cor} Let $g$ be a selfdual Einstein metric of nonzero scalar
curvature. Then there is a linear isomorphism from space of Killing fields of
$g$ to the space of twistors.
\end{cor}
This is not true if $g$ has zero scalar curvature, when $D^g$ is flat on
$\so\low_\asd(TM)$ and $g$ is locally hyper\kahl/. Then Proposition \ref{sdeK}
(i) says that $(D^gK)^\asd$ is parallel, and so $K$ is either triholomorphic,
or $(D^gK)^\asd$ is a nonzero constant multiple of one of the complex
structures.  Proposition \ref{sdeK} (ii) says that the Killing field
associated to a twistor is triholomorphic.

Let us remark, however, that this isomorphism does generalize to
quaternion-\kahl/ manifolds of nonzero scalar curvature~\cite{Sal:qkm},
and underlies the quaternion-\kahl/ quotient~\cite{GaLa:qro}.

\section{Surface-orthogonality}\label{so}

In this section we show that the orbits of two commuting Killing fields of a
selfdual Einstein metric with nonzero scalar curvature are necessarily
surface-orthogonal.

Recall that a Killing vector $K$ is hypersurface-orthogonal if and only if
$(D^gK)(X,Y)=0$ for all $X,Y$ orthogonal to $K$. Here we view $D^gK$ as a
$2$-form using the metric; $(D^gK)(X,Y)=-\frac12\ip{K,[X,Y]}$ and so
hypersurface-orthogonality means precisely that the orthogonal distribution is
integrable. Equivalently $K$ is hypersurface-orthogonal if and only if
$(*D^gK)(K)$ is zero. In four dimensions, this is a $1$-form called the
\emphdef{twist} of $K$.

Similarly, two linearly independent Killing vector fields $K,\tilde K$ are
(codimension $2$-) surface-orthogonal if and only if $(D^gK)(X,Y)=0$ and
$(D^g\tilde K)(X,Y)=0$ for all $X,Y$ orthogonal to both $K$ and $\tilde
K$. This means that the orthogonal distribution is integrable, and holds if
and only if $(*D^gK)(K,\tilde K)$ and $(*D^g\tilde K)(K,\tilde K)$ are both
zero. In four dimensions, these are both scalars, the \emphdef{twist scalars}.

We first collect some simple facts about commuting Killing fields.

\begin{lemma} Suppose $K, \tilde K$ are commuting Killing fields of a
$4$-dimensional metric $g$ and let $\Psi = (D^gK)^\asd$, $\tilde
\Psi=(D^g\tilde K)^\asd$. Then:
\begin{enumerate}
\item $(D^gK)(K,\tilde K)=0$ and $(D^g\tilde K)(K,\tilde K)=0$.
\item $d(|\tilde\Psi|^2)(K)=0$ and $d(|\Psi|^2)(\tilde K)=0$.
\end{enumerate}
\proofof{lemma} (i) Since $[K,\tilde K]=0$, we have $(D^g K)(\tilde K,\cdot)=
(D^g\tilde K)(K,\cdot)$. The results follow by contracting with $K$ and
$\tilde K$ respectively.

\noindent (ii) Since $\cL_Kg=0$ and $\cL_K\tilde K=0$, we have $\cL_KD^g\tilde
K=0$ and hence $\cL_K(\tilde\Psi)=0$. We deduce that $d(|\tilde\Psi|^2)
(K)=0$, and, in the same way, $d(|\Psi|^2)(\tilde K)=0$.
\end{proof}

Combining this with Proposition~\ref{sdeK} yields the surface-orthogonality.

\begin{prop}\label{surforth} Let $g$ be a selfdual Einstein metric of nonzero
scalar curvature, and suppose $K$, $\tilde K$ are linearly independent
commuting Killing fields.  Then the orthogonal distribution to
$\mathopen<\{K,\tilde K\}\mathclose>$ is integrable.
\proofof{prop} Contracting the formula of Proposition~\ref{sdeK} (i) with
$\Psi$, we obtain $\grad_g|\Psi|^2 = s^g\, \Psi K$. Now
$d(|\Psi|^2)(\tilde K)=0$ and $s^g$ is nonzero, so we deduce that
$\ip{\Psi(K),\tilde K}=0$. This implies that $(*D^gK)(K,\tilde
K)=(D^gK)(K,\tilde K)$ which vanishes by the above lemma.  A similar argument
shows that the other twist scalar also vanishes.
\end{proof}

\noindent In the zero scalar curvature case, the twist scalars are constant,
but they need not vanish.

\section{Scalar-flat \kahl/ metrics and Toda structures}\label{sfktoda}

On any Riemannian $4$-manifold $(M,g)$, the twistor equation has a geometric
interpretation due to Pontecorvo.

\begin{prop}\tcite{Pon:ths} Suppose that $\Psi$ is a section of
$\so\low_-(TM)$ satisfying the equation $D^g_X\Psi=(\gam\skwend X)^\asd$ for
all vector fields $X$, where $\gam$ is a $1$-form, and write, on the open set
where $\Psi$ is nonzero, $\Psi=f J$ where $J^2=-1$.

Then $(f^{-2}g,J)$ is a \textup(negatively oriented\textup) \kahl/ metric. In
particular the antiselfdual almost complex structure $J$ is integrable.
\proofof{prop} Contracting the equation with $\Psi$, we deduce that
$\gam=-2Jdf$. Hence $f D^g_XJ+df(X)J = -2(J df\skwend X)^\asd$, which may be
rewritten $D^g_XJ+(f^{-1}df\skwend X)\circ J- J\circ (f^{-1}df\skwend
X)=0$.  This means that $J$ is parallel with respect to the Levi-Civita
connection of $f^{-2}g$.
\end{proof}
Since the \kahl/ form of a negatively oriented \kahl/ metric $\hat g$ is
parallel, it is a twistor with respect to $\hat g$.  However, the twistor
equation is conformally invariant if $\Psi$ has weight $1$ (the other
component of the covariant derivative of $\Psi$, the divergence, is
equivalently the exterior derivative of the associated $2$-form, and so it is
conformally invariant if $\Psi$ has weight $-2$). This means that all
compatible \kahl/ metrics arise from twistors.

If $(M,g)$ is selfdual, then $(f^{-2}g,J)$ is a \emph{scalar-flat} \kahl/
metric~\cite{Gau:skc}.  Hence on a selfdual manifold, compatible scalar-flat
\kahl/ metrics are determined locally by solutions of a linear differential
equation.

Now suppose that $K$ is a Killing field of $g$, and that $\cL_K\Psi=0$.
Then $K$ is a holomorphic Killing field of the scalar-flat \kahl/ metric
$(f^{-2}g,J)$. LeBrun~\cite{LeBr:cp2} shows that such a scalar-flat \kahl/
metric $g_J$ is locally of the form
\begin{equation*}
g_J=w e^u(dx^2+dy^2)+w\, dz^2+w^{-1}(dt+A)^2,
\end{equation*}
where $\del_t$ is the Killing field, $u$ is a solution of the $\SU(\infty)$
Toda equation $u_{xx}+u_{yy}+(e^u)_{zz}=0$, and $w$ is a solution of its
linearization $w_{xx}+w_{yy}+(e^uw)_{zz}=0$, which is the compatibility
condition for the local existence of $A$ with $dA=w_x dy\wedge dz-w_y dx\wedge
dz +(e^uw)_z dx\wedge dy$. The scalar-flat \kahl/ metric is hyper\kahl/ if and
only if $u_z$ is a multiple of $w$, when LeBrun's construction reduces to that
of Boyer and Finlay~\cite{BoFi:kve} (or the Gibbons--Hawking Ansatz if
$u_z=0$).

A geometrical interpretation of this construction is obtained by relating it
to the Jones--Tod correspondence~\cite{JoTo:mew}. Given a selfdual space $M$
with a nonvanishing conformal vector field $K$, the local quotient of $M$ by
$K$ is a $3$-dimensional Einstein--Weyl space $B$---recall that this is a
conformal manifold equipped with a torsion-free conformal connection $D$ (a
\emphdef{Weyl connection}) such that the symmetric trace-free part of the
Ricci curvature of $D$ vanishes~\cite{CaPe:ewg}. Weyl connections on a
conformal manifold form an affine space modelled on the space of $1$-forms. In
the Jones--Tod construction, there is a unique compatible metric for which $K$
is a vector field of constant length, and $D$ differs from the Levi-Civita
connection of the quotient metric by a multiple of the twist of $K$ (which
descends to a $1$-form on $B$).

Conversely, given an Einstein--Weyl space $B$, selfdual spaces $M$ with a
conformal vector field fibering over $B$ are locally determined by solutions
of the abelian monopole equation $*Dw=dA$ (where $w$ is a section of $L^{-1}$,
i.e., a scalar of weight $-1$, and $A$ is a $1$-form).

In explicit terms, a Weyl structure may be specified by a choice of
representative metric $h$ and a $1$-form $\omega$ such that $Dh=-2\omega\tens
h$.  The Einstein--Weyl structure in LeBrun's construction is given by
\begin{equation}\label{toda}\begin{split}
h&=e^u(dx^2+dy^2)+dz^2\\
\omega&=-u_z dz.
\end{split}\end{equation}
An Einstein--Weyl space which can be written in this form, for some solution
$u$ of the $\SU(\infty)$ Toda equation, is said to be \emphdef{Toda}.
These Einstein--Weyl spaces were introduced by Ward~\cite{Ward:sut}.

The reduction to three dimensions of the twistor equation for $\Psi$
governing scalar-flat K\"ahler metrics gives a linear description of
compatible `Toda structures' on an Einstein--Weyl space~\cite{Cal:gte}.

\begin{defn} A \emphdef{Toda structure} on an Einstein--Weyl
is a section $\cX$ of $\smash{L^{-1/2}}\tens TB$ such that
$D\cX=\sigma\,\iden$ for some section $\sigma$ of $\smash{L^{-1/2}}$.  In
other words, $\cX$ is a weighted vector field with tracelike covariant
derivative.
\end{defn}
Using work of Tod~\cite{Tod:p3} (see also~\cite{CaPe:sdc}), it was shown
in~\cite{Cal:gte} that if $\cX$ is a nonvanishing Toda structure, then in the
gauge $(h,\omega)$ with $|\cX|=1$, called the \emphdef{LeBrun--Ward gauge},
the Einstein--Weyl space with is of the form~\eqref{toda}, for some solution
of the $\SU(\infty)$ Toda equation, and $\cX=\del_z$. 

We end this section by stating the crucial result of Tod~\cite{Tod:sut} which
characterizes selfdual Einstein metrics $g$ of nonzero scalar curvature with a
Killing field $K$. Proposition~\ref{sdeK} shows that these admit a solution
$\Psi$ of the twistor equation, and it is clear that $\cL_K\Psi=0$.

\begin{prop}\tcite{Tod:sut}\label{Todprop}
Let $g$ be a selfdual Einstein metric of nonzero scalar curvature with a
Killing field $K$. Then $g$ is locally isometric to $z^{-2}g_J$, where $g_J$
is a scalar-flat K\"ahler metric arising from LeBrun's construction with
$w=2-zu_z$. Conversely, on any Toda Einstein--Weyl space, $2-zu_z$ is a
solution of the abelian monopole equation and $z^{-2}g_J$ is Einstein.
\end{prop}

Note that a Toda structure only determines $z$ up to translation, so this
construction gives a one parameter family of selfdual Einstein metrics over
any Toda Einstein--Weyl space. However, it is very difficult to obtain
explicit solutions of the $\SU(\infty)$ Toda equation! In this paper, we are,
in effect, exploiting some implicit solutions of the $\SU(\infty)$ Toda
equation found by Ward~\cite{Ward:sut}.

\section{The Joyce Ansatz and Einstein--Weyl spaces}\label{jaew}

In~\cite{Joy:esd}, Joyce studied selfdual spaces with a surface-orthogonal
action of the $2$-torus by conformal transformations and constructed selfdual
conformal metrics on connected sums of complex projective planes.  To do this,
he first considered the local problem, and showed how selfdual conformal
metrics with a pair of surface-orthogonal commuting conformal vector fields
are generically determined by two solutions of a linear equation for a spinor
field $\Phi$ on the hyperbolic plane $\cH^2$.

On the other hand, in~\cite{Ward:sut}, Ward gave examples of Toda
Einstein--Weyl spaces by taking quotients of Gibbons--Hawking metrics
constructed from axisymmetric harmonic functions (AHFs) on $\R^3$. These
spaces were studied further in~\cite{Cal:gte}, where it was shown that they
are determined by a single solution $\Phi$ of the same linear equation on
$\cH^2$.

We refer to this equation for $\Phi$, which can be written
$\dbar\Phi=\frac12\overline\Phi$, as the \emphdef{Joyce equation}.

The two constructions can be related using the Jones--Tod correspondence.  To
do this in a symmetrical manner, we view the span of two commuting conformal
vector fields as a $2$-dimensional linear family $K^s$, $s\in\V$, where $\V$
is a $2$-dimensional real vector space on which we fix an area form (i.e.,
$\Wedge^2\V=\R$) once and for all. We say that $K^s$ as a \emphdef{pencil} of
conformal vector fields, since the nonzero elements up to scale are
parameterized by a real projective line.

We shall call a selfdual conformal manifold with a surface-orthogonal pencil
of conformal vector fields a \emphdef{Joyce space}. A pencil of solutions
$\boldsymbol\Phi=(\Phi_s)$ of the Joyce equation determines a Joyce space. On
the other hand, a single solution $\Phi$ of the Joyce equation defines an
Einstein--Weyl space \emphdef{with an axial symmetry}, i.e., admitting a
surface-orthogonal divergence-free conformal vector field $K$ preserving the
Weyl connection. We shall see that the pencil of quotients of a Joyce space
are the Einstein--Weyl spaces with an axial symmetry defined by the components
$\Phi_s$ of $\boldsymbol\Phi$ (in fact this pencil is parameterized by
$\V\dual$, not $\V$, which is why it is convenient to fix an area form on
$\V$). First let us summarize the two constructions.

\begin{prop} Suppose that $(N,g\low_N)$ is a hyperbolic $2$-manifold
with a spinor bundle $\cW$, i.e., a real rank $2$ vector bundle with a complex
structure such that $\cW\tens_\C\cW=TN$. Let $g\low_\cW$ be the induced
Hermitian metric on $\cW$.

\smallskip\noindent\textup{(i)~\cite{Cal:gte}} Suppose that
$\Phi\in\Cinf(N,\cW^*)$ is a solution of the Joyce equation which is
nonvanishing on an open subset $U$ of $N$, and let $\pi\colon B\to U$ be a
flat principal $S^1$- or $\R$-bundle with fibre coordinate $\psi$. Then
\begin{align*}
g&=|\Phi|^2\pi^*g\low_{N}+d\psi^2\\
\omega&=\Phi^2/|\Phi|^2
\end{align*}
is an Einstein--Weyl space with an axial symmetry. Conversely any connected
Einstein--Weyl space with an axial symmetry is either flat with translational
symmetry, or is locally isomorphic to one of these.

\smallskip\noindent\textup{(ii)~\cite{Joy:esd}} Suppose that
$\boldsymbol\Phi\in\Cinf(N,\cW^*)\tens\V$ is a pencil of solutions of the
Joyce equation which induces a positive isomorphism $\cW_x\to\V$ of real
vector spaces for all $x$ in an open subset $U$ of $N$, and let $\pi\colon
M\to U$ be a flat principal $\V/\Lambda$-bundle where $\Lambda$ is a discrete
subgroup of $(\V,+)$. Then the conformal class of the metric
\begin{equation*}
\pi^*g\low_N+g\low_\cW(\boldsymbol\Phi^{-1}(\cdot),\boldsymbol\Phi^{-1}(\cdot))
\end{equation*}
\textup(where we identify $TM$ with $\pi^*TN\dsum(M\cross\V)$, using the
principal connection\textup) is a Joyce space. Conversely, any connected Joyce
space is either locally conformally hyper\kahl/ with a pencil of
triholomorphic Killing fields, or is locally isomorphic to one of these.

\smallskip\noindent
\textup(For local questions we may as well take $N=\cH^2$, $B=U\cross\R$
and $M=U\cross\V$.\textup)
\end{prop}

The conventions used to identify $\dbar\Phi$ with $\frac12\overline\Phi$ in
the Joyce equation are crucial. If the curvature of the hyperbolic metric is
$-1$, the isomorphism from $\overline\cW{}^*$ (which is just $\cW^*$ with the
opposite complex structure) to $\overline{T\dual N}\tens_\C\cW^*$ must have
norm $1$ (as an isomorphism of real vector bundles) with respect to the
Hermitian metrics $g\low_\cW$ and $g\low_N$.

To clarify this, we follow Joyce by giving an explicit and purely real
interpretation of the Joyce equation. We take $N$ to be the hyperbolic plane
$\cH^2$ and introduce half-space coordinates $(\rho>0,\eta)$, so that
$g_{\cH^2}=(d\rho^2+d\eta^2)/\rho^2$. The metric on $\cW$ then has the form
$g\low_\cW=\mu_0^2+\mu_1^2$ where $\mu_0^2-\mu_1^2=d\rho/\rho$ and
$2\mu_0\mu_1=d\eta/\rho$ is the identification of $S^2_0\cW^*$ with $T\dual
N$.

We write the Joyce equation for the components $\Phi=A_0\mu_0+A_1\mu_1$
with respect to this orthonormal frame. The Levi-Civita connection of the
hyperbolic metric induces a Hermitian connection $\nabla$ on $\cW$ and it is
straightforward to compute the connection coefficients
$\nabla\mu_0=-\frac{d\eta}{2\rho}\tens \mu_1$ and
$\nabla\mu_1=\frac{d\eta}{2\rho}\tens \mu_0$. Hence
\begin{equation*}
\dbar\Phi=\bigl(\rho (A_0)_\rho+\rho (A_1)_\eta-\tfrac12A_0\bigr)\mu_0
-\bigl(\rho (A_1)_\rho-\rho (A_0)_\eta-\tfrac12A_1\bigr)\mu_1,
\end{equation*}
which equals $\frac12(A_0\mu_0+A_1\mu_1)$ if and only if
\begin{align*}
(A_0)_\rho+(A_1)_\eta&=A_0/\rho\\
(A_0)_\eta-(A_1)_\rho&=0.
\end{align*}

By solving one of these equations, we can reduce the Joyce equation to an
equation for a single function. The most obvious way to do this is to use the
second equation to set $A_0=G_\rho$ and $A_1=G_\eta$ so that the first
equation becomes $G_{\rho\rho}+G_{\eta\eta}=G_\rho/\rho$.  If we set
$G=\rho^{1/2}F$ then $F$ is an eigenfunction of the Laplacian on $\cH^2$ with
eigenvalue $3/4$.

Alternatively we can use the first equation to put $A_0=-\rho V_\eta$ and
$A_1=\rho V_\rho$, so that the second equation becomes $\rho
V_{\eta\eta}+(\rho V_\rho)_\rho=0$. This means that $V$ is an AHF on $\R^3$
with metric $d\rho^2+d\eta^2+\rho^2d\theta^2$.

By construction, the equation for $F$ is the integrability condition for $V$
and vice-versa, so the relation between $F$ and $V$ is a simple example of a
B\"acklund transformation.

Suppose now that we have two solutions $A_0\mu_0+A_1\mu_1$ and
$B_0\mu_0+B_1\mu_1$ of the Joyce equation. The corresponding Joyce space
has a compatible metric
\begin{equation}\label{eq:gen-metric}
g_0=(A_0 B_1 - A_1 B_0) g\low_{\cH^2}+
\frac{(A_0 d\phi - B_0 d\psi)^2+(A_1 d\phi - B_1 d\psi)^2}
{A_0 B_1 - A_1 B_0}.
\end{equation}
Via the Jones--Tod correspondence~\cite{JoTo:mew}, the quotient by $\del_\phi$
is an Einstein--Weyl space, and we can compute it by rediagonalizing $g_0$ and
rescaling by $(A_0^2+A_1^2)/(A_0 B_1 - A_1 B_0)$ to give
\begin{equation*}
(A_0^2+A_1^2)g\low_{\cH^2}+d\psi^2
+\left(\frac{A_0^2+A_1^2}{A_0 B_1 - A_1 B_0}\right)^2
\left(d\phi-\frac{(A_0 B_0 + A_1 B_1) d\psi} {A_0^2+A_1^2}\right)^2.
\end{equation*}
In this form, we can read off the Einstein--Weyl space
$(g\low_B,\omega\low_B)$ and abelian monopole $(\gmw,A)$, using the abelian
monopole equation $*(d\gmw-\omega\low_B\gmw)=dA$ to compute $\omega\low_B$.
The result is:
\begin{align*}
g\low_B&=(A_0^2+A_1^2)g\low_{\cH^2}+d\psi^2,&
\omega\low_B&=\frac{2A_0A_1\,d\eta+(A_0^2-A_1^2)d\rho}
{\rho(A_0^2+A_1^2)},\\
\gmw&=\frac{A_1 B_0 - A_0 B_1}{A_0^2 + A_1^2},&
A&=-\frac{A_0 B_0 + A_1 B_1} {A_0^2+A_1^2} d\psi.
\end{align*}
Hence, as expected, the result is the Einstein--Weyl space constructed from
the solution $A_0\mu_0+A_1\mu_1$ of the Joyce equation. However, we also see
explicitly how the second solution $B_0\mu_0+B_1\mu_1$ of the Joyce equation
determines an abelian monopole on this Einstein--Weyl space. In the remainder
of this section we discuss the scalar-flat K\"ahler metrics and Toda
structures that will enable us to characterize the case that $(M,g_0)$ is
conformally Einstein.

An important point is that \emph{any} Joyce space admits a family of
scalar-flat K\"ahler metrics~\cite{Joy:esd}, and \emph{any} Einstein--Weyl
space with an axial symmetry admits a family of Toda
structures~\cite{Cal:gte}. In fact, Joyce observes that each point at infinity
of $\cH^2$ determines a scalar-flat \kahl/ metric in the conformal class on
$M$: more precisely, for any half-space coordinates $(\rho,\eta)$ on $\cH^2$,
the metric $\rho g_0$ is scalar-flat K\"ahler.  By Pontecorvo's
work~\cite{Pon:ths}, this means that $M$ has a $2$-dimensional linear family
of solutions of the twistor equation, which we parameterize by a
$2$-dimensional real vector space $\W$.

Since these scalar-flat \kahl/ metrics are invariant under the entire pencil
of conformal vector fields, LeBrun's work~\cite{LeBr:cp2} shows that each
scalar-flat \kahl/ metric determines a Toda structure on each quotient
Einstein--Weyl space. This fits together with a further characterization of
Einstein--Weyl spaces with an axial symmetry~\cite{Cal:gte}: they are
precisely the Einstein--Weyl spaces admitting a pencil of Toda
structures. These Toda structures are parameterized by the same
$2$-dimensional real vector space $\W$, and are obtained from a choice of
half-space coordinates $(\rho,\eta)$ on $\cH^2$ as follows.

First introduce the functions $G,V$ with $G_\rho=A_0=-\rho V_\eta$ and
$G_\eta=A_1=\rho V_\rho$. (Note that $G$ is determined by $V$ up to
translation: one way to define $G$ is to choose an AHF $U$ with $U_\eta=V$, so
that $G=\rho U_\rho$.) Then, after a conformal rescaling by $\rho^2$, the
Einstein--Weyl structure may be written in the form~\eqref{toda},
\begin{equation*}
h=\rho^2(dV^2+d\psi^2)+dG^2,\qquad
\omega=\frac{2V_\eta}{\rho^2(V_\rho^2+V_\eta^2)} dG,
\end{equation*}
where we set $x=V$, $y=\psi$, $z=G$, and $e^u=\rho^2$.

We know that $u_z$ is an abelian monopole giving rise to hyper\kahl/
metric~\cite{BoFi:kve}, while the Einstein metrics we seek are obtained from a
$2-zu_z$ monopole (once we have fixed the Toda structure and the translational
freedom in the $z$ coordinate). The conformal rescaling by $\rho^2$ multiplies
the abelian monopoles $\gmw$ by (a constant multiple of) $1/\rho$, so in the
original gauge, these monopoles may be written
\begin{align*}
\frac\rho2 u_z&= \frac{A_0}{A_0^2+A_1^2}\\
\frac\rho2(2-zu_z)&=\frac{\rho(A_0^2+A_1^2)-G A_0}{A_0^2+A_1^2}.
\end{align*}
Clearly the first of these is the abelian monopole associated to the solution
$\mu_1$ of the Joyce equation. The second turns out to be associated to the
solution $B_0\mu_0+B_1\mu_1$ with $B_0=\rho A_1-\eta A_0$ and $B_1=G-\rho
A_0-\eta A_1$. This latter formula is rather mysterious at present, though one
easily checks that it does satisfy the Joyce equation, and that
$A_1B_0-A_0B_1=\rho(A_0^2+A_1^2)-G A_0$. The origin of this transformation of
the Joyce equation will be explained later.

\medbreak

\ifaddpics
\begin{center}
\includegraphics[width=.4\textwidth]{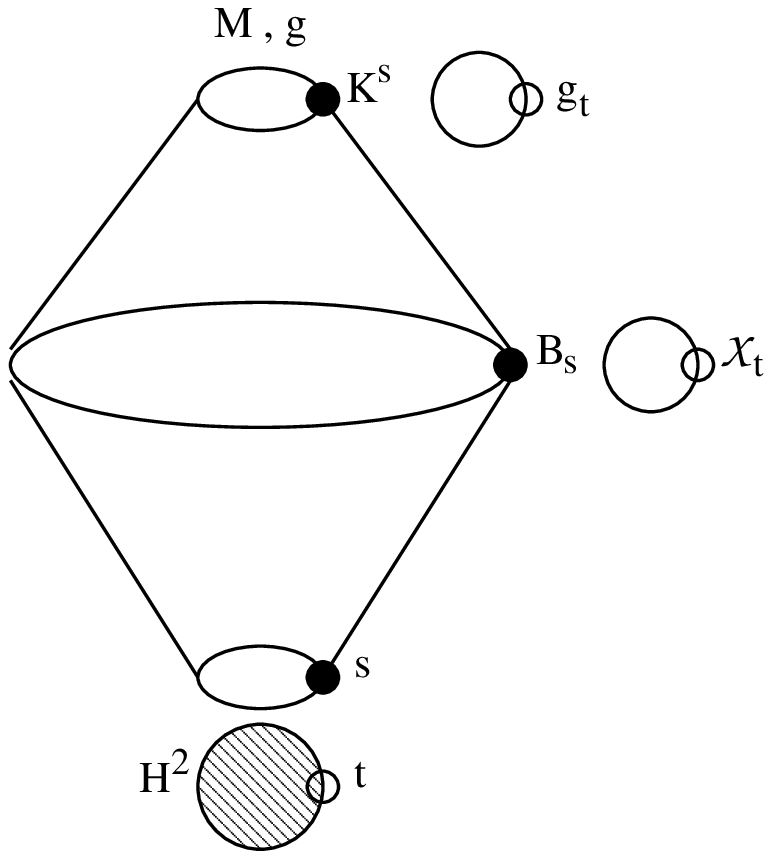}

\medbreak
Figure 1.
\end{center}
\fi

\bigbreak

Let us now summarize the discussion of this section (see Figure 1). First,
each solution $\Phi$ of the Joyce equation determines an Einstein--Weyl space
with an axial symmetry. A pencil of solutions $\boldsymbol\Phi=(\Phi_s)$
defines a pencil of Einstein--Weyl spaces $B_s$ and a selfdual space $M$ with a
surface-orthogonal pencil of conformal vector fields $K^s$, having the
Einstein--Weyl spaces $B_s$ as quotients. On the other hand, any such selfdual
space $M$ has a pencil of compatible scalar-flat \kahl/ metric $g_t$ and these
determine a pencil of Toda structures $\cX_t$ on each quotient Einstein--Weyl
space $B_s$. The choice of a scalar-flat \kahl/ metric or Toda structure is
given, up to homothety, by a point $t$ at infinity on the hyperbolic plane
$\cH^2$. We can fix the homothety freedom by introducing half-space
coordinates for the given point at infinity.  We shall see in the next section
how to regard points at infinity of $\cH^2$ as a pencil.

\section{Spin geometry of the hyperbolic plane}\label{sghp}

The hyperbolic plane $\cH^2$ is most naturally described as the space of
timelike lines in a $3$-dimensional Lorentzian vector space: it inherits a
Riemannian metric from the observation that each timelike line meets the
hyperboloid of two sheets at two points, one in each sheet; thus we may
identify $\cH^2$ with one of these sheets, with the induced metric.  Planar
models (such as the Poincar{\'e} disc or half-plane) are obtained by
introducing coordinates (pairs of real-valued functions with independent
differentials) on $\cH^2$.

In this section we describe the geometry of spinors on the hyperbolic plane by
equipping the Lorentzian vector space with some extra structure, following
Rees~\cite{Ree:nog} and Iversen~\cite{Ive:hg}.

Let $\W$ be a $2$-dimensional real vector space with an area form $\eps$ (so
$\Wedge^2\W=\R$). Then $S^2\W$ is a $(++-)$ Lorentzian vector space with
metric
$\ip{v_1v_2,w_1w_2}=-\eps(v_1,w_1)\eps(v_2,w_2)-\eps(v_1,w_2)\eps(v_2,w_1)$.
Alternatively, we can use $\eps$ to identify $S^2\W$ with $\sgl(\W)$, the Lie
algebra of traceless linear maps $A\colon\W\to\W$, so that the Lorentzian
quadratic form is $A\mapsto-\det A$, i.e., the timelike vectors have positive
determinant.

This linear algebra is a manifestation of the well-known isomorphism between
$\Spin(2,1)$ and $\SL(2,\R)$: $\W$ is the space of spinors of $S^2\W$.  The
hyperboloid of two sheets is the surface $\det A=1$ and so we use this
description to identify $\cH^2$ with the positive definite elements of
$S^2\W$ of determinant one.

In more concrete terms, choosing a unimodular basis of $\W$, we can
parameterize $\cH^2$ in $S^2\W$ by the matrices
\begin{equation*}
A(\rho,\eta)=\frac1\rho\left[\begin{matrix}1&\eta\\
\eta&\rho^2+\eta^2\end{matrix}\right].
\end{equation*}
Note that the corresponding traceless matrices are obtained by multiplying
these symmetric matrices by $J=\bigl[\begin{smallmatrix} 0 &-1\\ 1&
0\end{smallmatrix}\bigr]$, which has determinant one.  We easily compute that
\begin{align*}
dA &= \frac{d\rho}\rho\left[\begin{matrix}-1/\rho&-\eta/\rho\\
-\eta/\rho&(\rho^2-\eta^2)/\rho\end{matrix}\right]+
\frac{d\eta}\rho\left[\begin{matrix}0&1\\
1&2\eta\end{matrix}\right]\\
\tag*{so that}
\ip{dA,dA}&=\frac{d\rho^2+d\eta^2}{\rho^2}.
\end{align*}
Thus $\rho$ and $\eta$ are functions on $\cH^2$ identifying it with the
standard half-space model of the hyperbolic plane. They are only defined
once we have chosen the unimodular basis of $\W$.

The advantage of this model of the hyperbolic plane is that spinors are easy
to handle. Indeed we can identify $\cH^2\cross\W$ with the spinor bundle by
noting that for each $A\in\cH^2$, $A^{-1}$ is, by definition, a positive
definite unimodular inner product on $\W$, and this equips $\cW=\cH^2\cross\W$
with a metric.  The induced complex structure evidently satisfies
$\cW^2=T\cH^2$, since $\cW^2_A$ consists of the symmetric elements in
$\W\tens\low_\R\W$ which are traceless with respect to the inner product
$A^{-1}$, i.e., orthogonal to $A$, and this is the tangent plane $T_A\cH^2$.

In terms of a half-space model, a frame for $\cW$ is given by the vectors
\begin{equation*}
m_0=\Bigl[\begin{matrix}0\\ \sqrt\rho\,\end{matrix}\Bigr],\qquad
m_1=\Bigl[\begin{matrix}1\big/\!\sqrt\rho\,\\
                \eta\big/\!\sqrt\rho\,\end{matrix}\Bigr]
\end{equation*}
with dual frame
\begin{equation*}
\mu_0=\Bigl[\begin{matrix}\;\;\eta\big/\!\sqrt\rho\,\\
                  -1\big/\!\sqrt\rho\,\end{matrix}\Bigr],\qquad
\mu_1=\Bigl[\begin{matrix}\sqrt\rho\,\\0\end{matrix}\Bigr].
\end{equation*}
One easily sees that these are orthogonal and of unit length with respect to
$A$, and that
\begin{equation*}
m_0^2-m_1^2=dA(\rho\del_\rho),\qquad 2m_0m_1=dA(\rho\del_\eta).
\end{equation*}

We refer to the ``constant'' sections of $\cW=\cH^2\cross\W$ as
\emphdef{twistors}---they are certainly not parallel with respect to the
induced Hermitian connection on $\cW$, since the hyperbolic metric is not
flat. We shall not need to discuss what equation they satisfy, since we have
an explicit description of them in terms of the orthonormal frame
$m_0,m_1$: they are the sections of the form
\begin{equation*}
\Bigl[\begin{matrix}a\\b\end{matrix}\Bigr]=
-\frac{a\eta-b}{\sqrt\rho}m_0+a\sqrt\rho\,m_1,
\end{equation*}
for $a,b\in\R$. In particular, the norm of $[\begin{smallmatrix} a\\b
\end{smallmatrix}]$ is $\sqrt{a^2\rho^2+(a\eta-b)^2}/\sqrt\rho$.

For a geometric interpretation, note that $P(\W)$ is the space of null
lines in $S^2\W$, which are the points at infinity of the hyperbolic plane. A
point at infinity identifies the hyperbolic plane with the half-space model.
The inverse square norm of a twistor gives a `$\rho$' coordinate for this
half-space model. Evidently, the `$\rho$' coordinate determines the twistor up
to a sign.

\section{A tale of two pencils}\label{ttp}

We have introduced two pencils, parameterized by the vector spaces
$\V$ and $\W$. The first vector space $\V$ parameterizes the pencils of
\begin{itemize}
\item conformal vector fields on the selfdual space $M$ with torus
symmetry;
\item Einstein--Weyl spaces with an axial symmetry arising as quotients
of $M$;
\item solutions of the Joyce equation determining these Einstein--Weyl spaces.
\end{itemize}
The second vector space $\W$ parameterizes the pencils induced by
(deprojectivized) points at infinity of the hyperbolic plane, i.e., the
pencils of
\begin{itemize}
\item compatible scalar-flat \kahl/ metrics on $M$;
\item Toda structures on the Einstein--Weyl spaces;
\item twistors on the hyperbolic plane.
\end{itemize}
The idea now is that when the selfdual space $M$ is given by a selfdual
Einstein metric $g$ of nonzero scalar curvature, these pencils are the same:
\begin{itemize}
\item each Killing field determines a scalar-flat metric;
\item each Einstein--Weyl space has one of its Toda structures distinguished;
\end{itemize}
and therefore there must be a field on $\cH^2$ which determines a linear map
from twistors to solutions of the Joyce equation, and the resulting linear
family of solutions is the pencil defining the underlying conformal structure
of the Einstein metric. In this section we prove that this field is $F$. The
heart of the argument is the following result.
\begin{prop}\label{FtoPhi}
Let $F$ be an eigenfunction of the Laplacian on the hyperbolic plane with
eigenvalue $3/4$ and let $\varphi$ be a twistor.  Then $\Phi=\frac12 F
\,\flat\varphi+dF\act\varphi$ is a solution of the Joyce equation, where
$dF\act\varphi$ denotes the natural pairing $T^*\cH^2\tens\cW\to\cW^*$.

\proofof{prop} Given a twistor $\varphi$, we can use the freedom in
the choice of half-space model to set $\varphi=m_0/\sqrt\rho$.
In these half-space coordinates $F$ satisfies the equation
\begin{equation*}
F_{\rho\rho}+F_{\eta\eta}=\frac{3F}{4\rho^2}.
\end{equation*}
Since $dF=\rho F_\rho (\mu_0^2-\mu_1^2)+2\rho F_\eta \mu_0\mu_1$ direct
calculation yields $\Phi=(\rho^{1/2}F)_\rho\mu_0+ (\rho^{1/2}F)_\eta \mu_1$,
so $\Phi$ satisfies the Joyce equation.
\end{proof}

\begin{proof}[Proof of Theorem~\textup{\ref{main}}] Suppose that $g$ is a
selfdual Einstein metric of nonzero scalar curvature on $M$ with two commuting
Killing fields.  Let us review what we have proven so far about $g$. Firstly,
by Proposition~\ref{surforth}, the Killing fields are surface orthogonal, and
therefore, by~\cite{Joy:esd}, the conformal class of $g$ is a Joyce space.
(If $g$ is conformally hyper\kahl/, it must be conformally flat, i.e.,
locally isometric to $S^4$ or $\cH^4$, but then the Killing fields of $g$
cannot all be triholomorphic with respect to the flat hyper\kahl/ metric.)

The quotient of $g$ by one of its Killing fields is an Einstein--Weyl space
$B$ with an abelian monopole $\gmw$~\cite{JoTo:mew}. The work of
Tod~\cite{Tod:sut} shows that the choice of Killing field determines a
compatible scalar-flat \kahl/ metric on $M$, a Toda structure on $B$, and a
coordinate $z$ on $B$ such that the monopole $\gmw$ is $2-zu_z$. On the other
hand, we showed in section~\ref{jaew} that $B$ is an Einstein--Weyl space with
an axial symmetry. Since the Toda structure is invariant under this symmetry,
it is one of the Toda structures determined by a point at infinity on
$\cH^2$. Introducing compatible half-space coordinates we may write the
solution of the Joyce equation corresponding to this Einstein--Weyl space as
$\Phi=A_0\mu_0+A_1\mu_1$ and then $z=G$ for some function $G$ on $\cH^2$ with
$G_\rho=A_0$ and $G_\eta=A_1$.

Now set $F=\rho^{-1/2}G$. Then $F_{\rho\rho}+F_{\eta\eta}=\frac34 F/\rho^2$
and $\Phi$ is obtained by applying $F$ to the twistor $m_0/\sqrt\rho$ as in
Proposition~\ref{FtoPhi}. On the other hand, applying $F$ to $(-\eta\,
m_0+\rho\,m_1)/\sqrt\rho$ yields the solution $\tilde\Phi=( \rho A_1-\eta
A_0)\mu_0+(G-\rho A_0-\eta A_1)\mu_1$, which is precisely the solution needed
to construct the $2-zu_z$ monopole on $B$.

Hence $F$ generates the pencil of solutions of the Joyce equation yielding the
underlying conformal structure of the selfdual Einstein metric $g$. The
distinguished scalar-flat \kahl/ metric is $\rho g_0$ and rescaling this by
$1/z^2$, according to Proposition~\ref{Todprop}, we recover the selfdual
Einstein metric $F^{-2} g_0$. The explicit formula~\eqref{metric} is obtained
from~\eqref{eq:gen-metric} by direct substitution.

As we remarked in the introduction, the reader who is not convinced that we
really have encoded the entire selfdual Einstein condition in the construction
can easily verify this directly. Such calculations amount to reproving Tod's
result~\cite{Tod:sut} in this special case.
\end{proof}

In the next section we shall be able to obtain a better understanding of
formula~\eqref{metric} after studying the Swann bundle. We will also indicate
there how to check directly that the metric is selfdual and Einstein.

\section{The Swann bundle}\label{swann}

The \emphdef{Swann bundle} $\cU(M)$ of a selfdual Einstein manifold $(M,g)$
with nonzero scalar curvature is defined to be the principal $\CO(3)$-bundle
of conformal frames of $\Wedge^2T\dual_\asd M$. In~\cite{Swa:hkqk}, Swann
showed how to define a canonical (pseudo-)hyperk\"ahler metric on a similar
bundle over any quaternion-\kahl/ manifold.

The hyper\kahl/ structure on $\cU(M)$ is obtained as follows.  The Levi-Civita
connection induces a principal $\CO(3)$ connection on $\pi\colon\cU(M)\to
M$. The horizontal bundle of $\cU(M)$ is isomorphic to $\pi^*TM$, which has
three tautological $2$-forms (determined by the frame of $\Wedge^2T\dual_\asd
M$ at each point of $\cU(M)$), whereas the vertical bundle of $\cU(M)$ is
isomorphic to $\cU(M)\cross\HQ$, since $\CO(3)\cong\HQ\mult/\{\pm1\}$, and
this has a standard triple of $2$-forms. Adding suitable multiples of the
horizontal and vertical components gives the three symplectic forms.

In order to describe this more explicitly, identify $\cU(M)$ locally with
$M\cross\CO(3)$ by choosing a frame $\Theta$ of $\Wedge^2T\dual_\asd M$. We
view $\Theta$ as a $2$-form on $M$ with values in $\Imag\HQ$. The connection
on $\Wedge^2T\dual_\asd M$ is given by an $\Imag\HQ$-valued $1$-form $\omega$
satisfying
\begin{equation*}
d\Theta-\omega\wedge\Theta+\Theta\wedge\omega=0
\end{equation*}
where $(\omega\wedge\Theta)(X,Y)=\omega(X)\Theta(Y)-\omega(Y)\Theta(X)$ is the
usual wedge product of quaternion-valued forms. The selfdual Einstein equation
is now
\begin{equation*}
d\omega-\omega\wedge\omega+s\Theta=0
\end{equation*}
where $s$ is a constant---up to a positive numerical
factor it is the scalar curvature of the selfdual Einstein metric $g$.

Passing to a double cover, we have an $\HQ$-valued coordinate $q$ given by the
projection $M\cross\HQ\mult\to\HQ\mult$. The hyper\kahl/ metric is then
\begin{equation*}
\tilde g=s |q|^2 g+|dq+q\omega|^2
\end{equation*}
with $\Imag\HQ$-valued \kahl/ form $\Omega=s\, q\Theta \overline
q+(dq+q\omega) \wedge(\overline{dq+q\omega})$. An easy computation gives
\begin{equation*}\notag
d\Omega=dq\wedge(s\Theta-\omega\wedge\omega+d\omega)\overline q
+q(s\Theta+d\omega-\omega\wedge\omega)\wedge d\overline q
+q(s\,d\Theta+\omega\wedge d\omega-d\omega\wedge\omega)\overline q
\end{equation*}
which vanishes if $g$ is a selfdual Einstein metric.

Let us turn now to the examples of Theorem~\ref{main}. In this
case a frame for $\Wedge^2_\asd T\dual M$ is given by
\begin{align*}
\Theta&=
\frac1{F^2}\Bigl(\bigl(-\tfrac14 F^2 +\rho^2(F_\rho^2+F_\eta^2)\bigr)
\frac{d\rho\wedge d\eta}{\rho^2}+\alpha\wedge\beta\Bigr)\iI\\&\quad
+\frac1{F^2}\Bigl(\bigl(\rho F_\rho+\iI\rho F_\eta\bigr)(\alpha-\iI \beta)-
\tfrac12 F\,(\alpha+\iI\beta)\Bigr)\wedge\frac{d\rho-\iI d\eta}{\rho}\jJ
\end{align*}
where $\iI,\jJ,\kK$ are the imaginary quaternions. A tedious computation
is rewarded by a remarkably simple formula for the connection $1$-form:
\begin{equation*}
\omega=\frac1F\Bigl(-\rho F_\eta \frac{d\rho}\rho+\bigl(\tfrac12 F+\rho F_\rho
\bigr)\frac{d\eta}\rho\Bigr)\iI - \frac1F (\alpha-\iI\beta)\jJ.
\end{equation*}
Computing $d\omega-\omega\wedge\omega$, we can check that the metric $g$ is
selfdual and Einstein with $s=1$.

Since the construction of the Swann bundle is canonical, the commuting Killing
fields of $g$ lift to give commuting trihamiltonian vector fields of $\tilde
g$. Now any hyperk\"ahler $8$-manifold with two commuting trihamiltonian
vector fields is given explicitly by a generalized Gibbons--Hawking
Ansatz~\cite{HKLR:hkms,PePo:hkm}: it is isometric to
\begin{equation*}
\Phi_{ij}\ip{d\boldsymbol x_i,d\boldsymbol x_j}
+{\Phi^{-1}}\low_{\!ij}(dt_i+A_i)(dt_j+A_j)
\end{equation*}
where $(\Phi_{ij},A_i)$ is a solution of a generalized abelian monopole
equation on $\R^2\tens\Imag \HQ$, whose coordinates $(\boldsymbol
x_1,\boldsymbol x_2)$ are the $\Imag\HQ$-valued momentum maps of the
trihamiltonian vector fields $(\del_{t_1},\del_{t_2})$.

In more invariant language, the matrix $\underline\Phi$ is a section of
$S^2\V$ over $\V\dual\tens\Imag\HQ$, where $\V$ is the $2$-dimensional real
vector space of Killing fields. For our metrics $\V\dual=\W$ and it is
convenient to use the non-constant frame dual to $(\alpha,\beta)$ for $\W$.

\begin{prop}\label{swannprop} Let $F$ be a local eigenfunction of the
Laplacian on $\cH^2$.  Then the hyper\kahl/ metric on the Swann bundle of the
associated selfdual Einstein metric is given by the generalized
Gibbons--Hawking Ansatz with
\begin{equation*}
\underline \Phi = \frac F{|q|^2}
\begin{pmatrix} \tfrac12 F+\rho F_\rho & \rho F_\eta\\
\rho F_\eta & \frac12 F-\rho F_\rho
\end{pmatrix}.
\end{equation*}
Furthermore, the $\Imag\HQ$-valued momentum maps of $\del_\psi$ and
$\del_\phi$ are
\begin{equation*}
\boldsymbol x_\psi=\frac{q\kK \overline q}{\sqrt\rho F}
\qquad\text{and}\qquad
\boldsymbol x_\phi=\frac{q(\eta+\rho\iI)\kK\overline q}{\sqrt\rho F}.
\end{equation*}
\proofof{prop} For the first part we show that $\underline\Phi^{-1}$ is
the metric on $\cU(M)\cross\W$ induced by $h$. Let $Q$ be the
matrix $\Bigl(\begin{smallmatrix} \frac12 F-\rho F_\rho & -\rho F_\eta\\
-\rho F_\eta & \frac12 F+\rho F_\rho
\end{smallmatrix}\Bigr)$, so that the selfdual Einstein metric is
\begin{equation*}
g=\frac1{F^2}\Bigl[\detm Q\, g\low_{\cH^2}
+\frac1{\detm Q}(\alpha\;\;\beta)
Q^2\Bigl(\begin{matrix}\alpha\\ \beta\end{matrix}\Bigr)\Bigr]
\end{equation*}
Then the metric on the torus induced by $\tilde g=|q|^2g+|dq+q\Theta|^2$ is
\begin{equation*}
\frac{|q|^2}{F^2\detm Q}(\alpha\;\;\beta)\bigl(\detm Q\,\iden+Q^2\bigr)
\Bigl(\begin{matrix}\alpha\\ \beta\end{matrix}\Bigr)
\end{equation*}
and $\detm Q\,\iden+Q^2=(\trace Q)\,Q=F\,Q$ by the Cayley--Hamilton theorem.
This is what we want, because $\underline\Phi^{-1}=|q|^2 Q/(F\detm Q)$.

For the momentum maps, we must compute the contraction of the
$\Imag\HQ$-valued symplectic form $q\Theta \overline
q+(dq+q\omega)\wedge(d\overline q-\omega\overline q)$ with $\del_\psi$ and
$\del_\phi$. The contraction with any vector field $X$ in the torus is
$q\bigl(\Theta(X)+\omega(X)\omega-\omega\,\omega(X)\bigr)\overline q
+q\omega(X)d\overline q+dq\, \omega(X)\overline q$.  It is straightforward to
compare this with $d\boldsymbol x_\psi$ or $d\boldsymbol x_\phi$ when
$X=\del_\psi$ or $\del_\phi$: in particular note that
$\omega(\del_\psi)=\kK/\sqrt\rho F$ and
$\omega(\del_\phi)=(\eta+\rho\iI)\kK/\sqrt\rho F$.
\end{proof}

We summarize this discussion of the Swann bundle by the following diagram.
\begin{diagram}[size=2em]
\cU(M)^8&\rTo^{\CO(3)}&M^4\\
\dTo^{T^2}&  &\dTo_{T^2}\\
\Imag\HQ\tens_0\W&\rTo^{\CO(3)} &\cH^2
\end{diagram}
Here we denote by $\Imag\HQ\tens_0\W$ the elements of $\Imag\HQ\tens\W$
which are not decomposable. In terms of a basis for $\W$, this means the
points where the coordinates $\boldsymbol x_1$ and $\boldsymbol x_2$ are
linearly independent.  Clearly the momentum maps $\boldsymbol x_\phi$
and $\boldsymbol x_\psi$ are linearly independent.

In order to understand the bottom arrow in this diagram, consider the
Grammian map $\Imag\HQ\tens\W\to S^2\W$ given in components by
\begin{equation*}
(\boldsymbol x_1,\boldsymbol x_2)\mapsto\left[\begin{matrix}
|\boldsymbol x_1|^2&\ip{\boldsymbol x_1,\boldsymbol x_2}\\
\ip{\boldsymbol x_1,\boldsymbol x_2}& |\boldsymbol x_2|^2
\end{matrix}\right].
\end{equation*}
The determinant of this matrix is $|\boldsymbol x_1\wedge\boldsymbol x_2|^2$,
which is nonnegative, and vanishes if and only if $\boldsymbol x_1$ and
$\boldsymbol x_2$ are linearly dependent. Otherwise the matrix is positive
definite, and so on dividing by $|\boldsymbol x_1\wedge\boldsymbol x_2|$, we
get a well-defined map $\Imag\HQ\tens_0\W\to \cH^2$.  This map is
$\CO(3)$-invariant, where $\CO(3)$ acts diagonally on
$\Imag\HQ\tens\W\cong\Imag\HQ\dsum\Imag\HQ$, and applying it to $(\boldsymbol
x_\phi,\boldsymbol x_\psi)$ gives the matrix
\begin{equation*}
A(\rho,\eta)=\frac1\rho\left[\begin{matrix}1&\eta\\
\eta&\rho^2+\eta^2\end{matrix}\right]
\end{equation*}
so the diagram commutes.

Finally in this section, we justify our use of the letter $\Phi$ both for
solutions of the Joyce equation defining selfdual Einstein metrics, and for
generalized monopoles defining their Swann bundles. Given any hyperbolic
eigenfunction $F$ on $\cH^2$, define $\tilde F\colon S^2\W_+\to\R$, where
$S^2\W_+$ denote the space of timelike elements of $S^2\W$ (i.e., matrices of
positive determinant), by requiring that $\tilde F$ has homogeneity $1/2$,
i.e., $\tilde F(\lambda v)=\lambda^{1/2}\tilde F(v)$, and that $\tilde
F\restr{\cH^2}=F$. Then $d\tilde F$ is a function on $S^2\W_+$ with values in
$S^2\W$ and its matrix with respect to the (homogeneity $1/2$) orthonormal
frame $(m_0,m_1)$ of $S^2\W_+\cross\W$ is
\begin{equation*}
\frac1{\sqrt{\detm A}}\begin{pmatrix} \tfrac12 F+\rho F_\rho & \rho F_\eta\\
\rho F_\eta & \frac12 F-\rho F_\rho
\end{pmatrix}\qquad\text{at}
\qquad A=\frac{\sqrt{\detm A}}{\rho}\left[\begin{matrix}1&\eta\\
\eta&\rho^2+\eta^2\end{matrix}\right]\in S^2\W_+.
\end{equation*}
\begin{itemize}
\item Pulling $d\tilde F$ back to $\Imag\HQ\tens_0\W$ gives the generalized
monopole $\underline\Phi$, since $\sqrt{\det A}$ pulls back to
$|\boldsymbol x_\psi\wedge \boldsymbol x_\phi|=|q|^2/F$.
\item Restricting $d\tilde F$ to $\cH^2$ gives the pencil of solutions
$\boldsymbol\Phi$ of the Joyce equation and this explains the form of
the metric in~\eqref{metric}.
\end{itemize}
Hence $\underline\Phi$ and $\boldsymbol\Phi$ are different manifestations of
the same object. Furthermore, this description in terms of a homogeneity $1/2$
function $\tilde F$ gives a natural interpretation of the equation for $F$:
$F=\tilde F\restr{\cH^2}$ is a hyperbolic eigenfunction with eigenvalue $3/4$
if and only if $\tilde F$ is a (homogeneity $1/2$) solution of the wave
equation.

\section{Quaternion-\kahl/ and hyper\kahl/ quotients}\label{qkhk}

Our motivation for constructing the Swann bundle of the selfdual Einstein
metrics of this paper is to provide an explicit relation between these metrics
and the quaternion-\kahl/ quotients of quaternionic projective space
$\HP{m-1}$. In~\cite{Gal:gmm,GaLa:qro} Galicki and Lawson defined an analogue
of the hyper\kahl/ quotient~\cite{HKLR:hkms} in quaternion-\kahl/ geometry, in
which the quotient of a $4(m-1)$-dimensional quaternion-\kahl/ manifold by a
$k$-dimensional group of symmetries is (at least locally) a
$4(m-k-1)$-dimensional quaternion-\kahl/ manifold.

This is of interest here, because quaternion-\kahl/ quotients of $\HP{m-1}$ by
an $(m-2)$-dimensional subtorus of a maximal torus $T^m$ in $\Symp(m)$ are
selfdual Einstein metrics of positive scalar curvature with $T^2$ symmetry.
These quotients were first studied by Galicki--Lawson~\cite{GaLa:qro} and
Boyer--Galicki--Mann--Rees~\cite{BGMR:3s7}; they are globally defined on
compact orbifolds. Hence if we could obtain explicitly the relation between
these quaternion-\kahl/ quotients and Theorem~\ref{main}, then we would have
explicit formulae for the (hithertoo only implicit) Galicki--Lawson metrics
and their generalizations.

Quaternion-\kahl/ quotients of a quaternion-\kahl/ manifold $Q$ may be related
to hyper\kahl/ quotients of its Swann bundle $\cU(Q)$. Indeed the symmetry
group lifts to an action on $\cU(Q)$ preserving the hyper\kahl/ structure, and
commuting with the $\CO(3)$-action.  The momentum map of this action is a
$\CO(3)$- and $G$-equivariant map $\mu\colon\cU(Q)\to\Imag\HQ\tens\lie{g}^*$,
where $\Imag\HQ$ carries the standard representation of $\CO(3)$ and $\lie{g}$
is the Lie algebra of $G$, so that $\lie{g}^*$ is the coadjoint
representation. The hyper\kahl/ quotient of $\cU(Q)$, given by
$\cU(Q)/\!/\!/G=\mu^{-1}(0)/G$, is therefore hyper\kahl/ with a
$\CO(3)$-action.  (It may or may not be a manifold, but in any case the
geometry of the local quotient is well-defined.)  Swann proved that
$\cU(Q)/\!/\!/G$ is the Swann bundle of the quaternion-\kahl/ quotient of
$Q$~\cite{Swa:hkqk}.  Indeed, when taking quotients of $\HP{m-1}$ one often
works in homogeneous coordinates, and this amounts to working on (the double
cover of) the Swann bundle $(\HQ^m\punc0)/\{\pm1\}$.

Now hyper\kahl/ quotients of $\HQ^m$ by tori are well understood.  Let $T^m$
be the maximal torus of $\Symp(m)$ acting on $\HQ^m$ by $(q_1,\ldots
q_m)\mapsto(e^{\iI t_1}q_1,\ldots e^{\iI t_m} q_m)$. We can describe an
$(m-2)$-subtorus of $T^m$ by declaring that its Lie algebra is the kernel of a
map $\R^m\to\R^2$ sending the standard basis $e_1,\ldots e_m$ (the generators
of the chosen $m$-torus) to some given $\alpha_1,\ldots\alpha_m\in\R^2$.
Evidently $\alpha_1,\ldots\alpha_m$ must be rational (up to an overall factor)
in order that the kernel of the map $\R^m\to\R^2$ is the Lie algebra of a
subtorus. However, even without this condition, we can still consider the
local hyper\kahl/ quotient of $\HQ^m$ by an $(m-2)$-dimensional family of
commuting triholomorphic Killing fields.

Specializing a result of Bielawski--Dancer~\cite{BiDa:gtt} to the case of
interest, we learn that the hyper\kahl/ quotient of $\HQ^m\punc0$ by this
subtorus (with zero momentum map in order to obtain the Swann bundle of a
quaternion-\kahl/ quotient) is given by the generalized Gibbons--Hawking
Ansatz with
\begin{equation*}
\Phi_{ij}=\sum_{k=1}^m \frac{(\alpha_k)_i(\alpha_k)_j}{r_k},
\end{equation*}
where $r_k = |(\alpha_k)_1\boldsymbol x_1+(\alpha_k)_2\boldsymbol x_2|$.

Note that $\R^2$ here is the Lie algebra of the quotient torus. In our setting
this is the vector space $\W$. Then $\alpha_1,\ldots\alpha_m$ define $m$
twistors on $\cH^2$, which we write as $\varphi_k= \bigl[\begin{smallmatrix}
a_k\\b_k \end{smallmatrix}\bigr]$, where $a_k=-(\alpha_k)_2$,
$b_k=(\alpha_k)_1$.  These $m$ twistors must determine the hyperbolic
eigenfunction $F$ in some way. Since $\underline\Phi$ is $d\tilde F$, we
compute that
\begin{equation*}
\tilde F=\sum_{k=1}^m r_k=\sum_{k=1}^m|b_k\boldsymbol x_1-a_k\boldsymbol x_2|
\end{equation*}
and therefore
\begin{equation*}
F=\sum_{k=1}^m \frac{\sqrt{\smash[b]{a_k^2}\rho^2
+(a_k\eta-b_k)^2}}{\sqrt\rho} = \sum_{k=1}^m |\varphi_k|
\end{equation*}
where $|\varphi_k|$ is the pointwise norm of the twistor
$\bigl[\begin{smallmatrix} a_k\\b_k \end{smallmatrix}\bigr]$. Indeed, it is
straightforward to check that the norm of a twistor $\varphi$ is a hyperbolic
eigenfunction: without loss of generality we can take
$\varphi=\bigl[\begin{smallmatrix} 0\\1 \end{smallmatrix}\bigr]
=m_0/\sqrt\rho$, so that $F=1/\sqrt\rho$.

For this `monopole solution' $F$, the pencil of solutions of the Joyce
equation degenerates: applying $F$ to
$\bigl(-(a\eta-b)m_0+a\rho m_1\bigr)/\sqrt\rho$, using
Proposition~\ref{FtoPhi}, gives $\Phi=2a\mu_1$, which is the solution of the
Joyce equation found by Joyce~\cite{Joy:esd}. Joyce superposed this solution
with its image under hyperbolic isometries in order to obtain his explicit
metrics on $n\CP2$.

The same trick yields interesting metrics here. For the monopole solution,
$F=1/\sqrt\rho$, $\frac14F^2-\rho^2(F_\rho^2 + F_\eta^2)$ is identically zero,
and so we only obtain a selfdual Einstein metric only on the empty set!
However for $m>1$ the `$m$-pole' solutions $F=\sum_{k=1}^m |\varphi_k|$ yield
quaternion-\kahl/ quotients of $\HP{m-1}$. Let us state our result more
precisely.

\begin{thm}\label{qkposthm} Let $\R^m$ be the Lie algebra of the maximal
torus of $\Symp(m)$ which acts on $\HP{m-1}$ by $[q_1\ratio\cdots\ratio
q_m]\mapsto[e^{\iI t_1}q_1\ratio\cdots\ratio e^{\iI t_m} q_m]$. Let $M^4$ be
the local quaternion-\kahl/ quotient of $\HP{m-1}$ by the $(m-2)$-dimensional
family of Killing fields in the kernel of the map $\R^m\to\W$ sending the
standard basis $e_1,\ldots e_m$ to $\varphi_1,\ldots\varphi_m\in\W$.

Then the selfdual Einstein metric on $M^4$ is given by~\textup{\eqref{metric}}
with $F=\sum_{k=1}^m|\varphi_k|$
\end{thm}
The solutions corresponding to reductions of $\HP{m-1}$ by an $(m-2)$-torus
yield the compact selfdual Einstein orbifolds we seek. These global reductions
arise when $\varphi_1,\ldots\varphi_m\in\W$ span a $2$-dimensional vector
space over the rationals, i.e., when they can be chosen to have rational
components.

The description of these metrics in terms of the hyperbolic plane links the
geometry to the topological analysis of
Boyer--Galicki--Mann--Rees~\cite{BGMR:3s7}, who describe a $T^2$-invariant
cell decomposition of their orbifolds over the closed disc, where the
principal orbits fibre over the open disc and the special orbits fibre over
the boundary. If we identify the open disc with the hyperbolic plane $\cH^2$,
then the boundary is the circle at infinity $P(\W)$. The $m$ twistors
determine $m$ marked points on this circle, corresponding to the fixed points
of the action.

Non-compact analogues of these metrics may be obtained by taking
quaternion-\kahl/ quotients of quaternionic hyperboloids $\HQ\cH^{p-1,q}$ and
$\HQ\cH^{p,q-1}$ by an $(m-2)$-subtorus of a Cartan subgroup of $\Symp(p,q)$,
$p+q=m$.  This includes, for instance, quotients of quaternionic hyperbolic
space $\HQ\cH^{m-1}$ when $\{p,q\}=\{m-1,1\}$~\cite{Gal:nmc,Gal:mcm}.

Such metrics may also be viewed as analytic continuations. Indeed, both the
constructions of this paper and the hyper\kahl/ and quaternion-\kahl/
quotients may be carried out in the holomorphic category, and from this point
of view real metrics are obtained by taking real slices of such complex
metrics.  In the holomorphic category, $|\varphi|$ must be taken to be a
choice of branch of the square root of the complex bilinear pointwise inner
product of $\varphi$ with itself. Other real slices are obtained by replacing
a sum of real twistors either by a sum of complex conjugate twistors, or by a
difference of real twistors, leading to multipole solutions of the form
\begin{equation*}
F=\sum_{k=1}^{r}
\bigl(|\varphi_{2k-1}+\iI\varphi_{2k}|+|\varphi_{2k-1}-\iI\varphi_{2k}|\bigr)
+\sum_{k=r+1}^{p}\bigl(|\varphi_{2k-1}|-|\varphi_{2k}|\bigr)
+\sum_{k=2p+1}^{p+q}|\varphi_{k}|,
\end{equation*}
where $\varphi_1,\ldots\varphi_{p+q}$ are real and $0\leq r\leq p$. (Note that
the norm of a complex twistor vanishes on the hyperbolic plane and so the
complex conjugate pairs are only defined on a branched cover.)

The freedom here corresponds to the fact that $\Symp(p,q)$ does not have a
unique Cartan subgroup up to conjugacy: for each $\Symp(1,1)$ factor, we can
take either an $S^1\cross\R$ subgroup or a $T^2$ subgroup, using (for
instances) elements of the form
\begin{equation*}
\begin{pmatrix} e^{\iI t_1}\cosh t_2&e^{\iI t_1}\sinh t_2\\
e^{\iI t_1}\sinh t_2&e^{\iI t_1}\cosh t_2 \end{pmatrix}
\qquad\text{or}\qquad
\begin{pmatrix} e^{\iI t_1}&0\\0&e^{\iI t_2} \end{pmatrix}
\end{equation*}
in $\Un(1,1)$. Hence if we suppose $p<q$ and write
$\HQ^{p,q}=(\HQ^{1,1})^p\cross\HQ^{q-p}$, then we obtain a Cartan subgroup of
$\Symp(1,1)\cross\cdots\cross\Symp(1,1)\cross\Symp(q-p)$ of the form
\begin{equation*}
(S^1\cross\R)^r\cross T^{2(p-r)}\cross T^{q-p}\leq
\Symp(1,1)^r\cross\Symp(1,1)^{p-r}\cross\Symp(q-p).
\end{equation*}
We should remark, however, that in addition to Cartan subgroups, there are
also maximal abelian subgroups containing nilpotent elements. We believe these
yield a mild generalization of $m$-pole solutions, in which dipoles can become
infinitesimal in a limiting process: $\lim_{\eps\to 0}
\frac1\eps(|\varphi_1+\eps\varphi_2|-|\varphi_1-\eps\varphi_2|)$.  We shall
not study this kind of solution here.

The global behaviour of the $m$-pole solutions can be approached either via
the quaternion-\kahl/ quotient, as in~\cite{BGMR:3s7,Gal:nmc,Gal:mcm}, or via
compactification arguments based on local models, as discussed by
Joyce~\cite{Joy:esd}. A detailed analysis of this would take us to far afield
here, so we turn instead to examples.

\section{Examples}\label{es}

In this section we study the simplest nontrivial examples of multipole
selfdual Einstein metrics. An $m$-pole solution $F$ is determined, up to sign
and reality choices, by $m$ distinct elements of the $2$-dimensional vector
space $\W$ of twistors, and conversely $F$ determines the twistors up to
sign. Concretely, with respect to a choice of unimodular basis, we may write
these elements as $\bigl[\begin{smallmatrix} a_k\\b_k
\end{smallmatrix}\bigr]$. Now the group $\SL(\W)$ acts naturally on the set
of $m$-tuples of elements of $\W$, yielding equivalent solutions, so there is
really only a $2m-3$ parameter family of solutions. Furthermore, the
solutions $F$ and $\lambda F$ yield the same Einstein metric for any
$\lambda\neq0$, so the moduli space actually has dimension $2m-4$.

We first remark that the dipole solutions ($m=2$) yield only hyperbolic and
spherical metrics. More precisely, using the $\SL(\W)$ freedom, we take
the two twistors to be $[\begin{smallmatrix} a\\0
\end{smallmatrix}]$ and $[\begin{smallmatrix} 0\\a \end{smallmatrix}]$.
Allowing for the homothety freedom, we are then left with three choices for
the solution:
\begin{align*}
F^+ &= \frac{1+\sqrt{\rho^2+\eta^2}}{\sqrt\rho},\qquad
F^- = \frac{1-\sqrt{\rho^2+\eta^2}}{\sqrt\rho}\\ \tag*{or}
F^c &= \frac{\sqrt{\rho^2+(\eta+\iI)^2}+\sqrt{\rho^2+(\eta-\iI)^2}}
{\sqrt\rho}.
\end{align*}
The first of these gives the spherical metric, while the other two give the
hyperbolic metric with inequivalent torus actions.

We now consider the case $m=3$, when the moduli space is $2$-dimensional.
This case is particular interesting since selfdual Einstein metrics in
this family have been studied in many places by diverse methods.
\begin{enumerate}
\item In the context of finding selfdual Einstein metrics of negative scalar
curvature with prescribed conformal infinity~\cite{LeBr:hcc}, the second
author found complete examples on the $4$-ball~\cite{Ped:emm}, depending on a
single parameter in $(-1,\infty)$ (denoted $m^2$ there). It was later
realized~\cite{LeBr:pac,Hit:tem} that when this parameter is $(2-n)/n$ (for
$n\in\Z$, $n\geq3$), analytic continuations of these metrics are complete on
$\cO(n)\to\CP1$ (and are conformally related to scalar-flat K\"ahler metrics
on $\cO(-n)$).
\item By taking quaternion-\kahl/ quotients of $\HP2$ by $S^1$, Galicki and
Lawson~\cite{GaLa:qro} found selfdual Einstein metrics of positive scalar
curvature on certain compact orbifolds $\cO_{q,p}$ ($p,q$ coprime with
$0<q/p\leq 1$). Negative scalar curvature analogues and generalizations have
also been studied~\cite{Gal:nmc,Gal:mcm}.
\item A subfamily of the $m=3$ metrics have local cohomogeneity one, and are
therefore bi-axial Bianchi metrics, which have been studied in many places, in
particular~\cite{BBe:nvr}. The metrics in (i) are all in this
subfamily. Furthermore, the quaternion-\kahl/ quotients considered in detail
in (ii) are mainly the local cohomogeneity one examples, although the general
case is analogous~\cite[Remark 4.27]{GaLa:qro}.
\item Apostolov and Gauduchon~\cite{ApGa:seh} classify explicitly selfdual
Einstein Hermitian metrics, i.e., admitting a selfdual complex structure: such
metrics automatically have torus symmetry. They are conformal to selfdual
K\"ahler metrics, which have been classified by Bryant~\cite{Bry:bkm} as the
specialization to four dimensions of K\"ahler metrics with vanishing Bochner
tensor. Conversely a generic selfdual K\"ahler metric is locally conformally
Einstein.  Apostolov and Gauduchon also show that quaternion-\kahl/ quotients
of $\HP2$ and quaternionic hyperbolic space $\HQ\cH^2$ by $S^1$ or $\R$ are
selfdual Hermitian (and the same is true for $\HQ\cH^{1,1}$). In particular
the $3$-pole solutions are all Hermitian.
\end{enumerate}
In view of all this work, we cannot claim that the $m=3$ examples are new.
Nevertheless, in addition to presenting a unified treatment, we are able to
give explicitly the parameter values yielding complete $3$-pole metrics on the
$4$-ball.  At the end of the section, we shall use a perturbation argument to
obtain complete $m$-pole metrics on the $4$-ball for any $m$, showing that the
moduli space is infinite dimensional.

Consider then the general $3$-pole solution. Using the $\SL(\W)$ freedom
we may write this in the form:
\begin{equation*}
F= \frac a{\sqrt\rho}
+\frac{b+c/m}2\frac{\sqrt{\rho^2+(\eta+m)^2}}{\sqrt\rho}
+\frac{b-c/m}2\frac{\sqrt{\rho^2+(\eta-m)^2}}{\sqrt\rho},
\end{equation*}
where $|m|=1$, but $m$ can be imaginary or real, $-m^2=\pm1$.  We refer to
these as Type I and Type II solutions respectively, after the Eguchi--Hanson I
and II metrics. It is convenient to work in (Eguchi--Hanson)-like coordinates
\begin{equation*}
\rho=\sqrt{R^2\pm1}\,\cos\theta,\qquad \eta=R\sin\theta
\end{equation*}
with $\theta\in(-\pi/2,\pi/2)$ so that
\begin{align*}
\sqrt\rho F &= a+b R+c\sin\theta\\ \tag*{and}
\rho^{-1}\bigl(\tfrac14 F^2-\rho^2(F_\rho^2 + F_\eta^2)\bigr)
&=\frac{b(aR\mp b)+c(a\sin\theta+c)}{R^2\pm\sin^2\theta}.
\end{align*}
Note that the zero-set of $F$ is a conformal infinity of the selfdual Einstein
metric, which has negative scalar curvature there. On the other hand the
zero-set of $\tfrac14 F^2-\rho^2(F_\rho^2 + F_\eta^2)$ is a singularity
separating domains of positive and negative scalar curvature, and the metric
is incomplete there. Let us consider the Type I and Type II metrics
separately.

In the Type I case, $F$ is only globally defined on a branched double cover of
the hyperbolic plane: we can regard $R=0$ be the branch cut between $R\geq 0$
and $R\leq 0$. Without loss of generality, we can assume $a$ is nonzero and
use the homothety freedom to set $a=1$. We can also suppose $b,c\geq0$ (using
$R\mapsto-R$ and $\theta\mapsto-\theta$).

When $b$ is nonzero we have, for each $\theta\in(-\pi/2,\pi/2)$, a unique
value of $R$, namely $R_\infty=-(1+c\sin\theta)/b$, at which $F=0$, and a
unique value $R_\pm=(b^2+c^2+c\sin\theta)/b$ where $\frac14
F^2-\rho^2(F_\rho^2 + F_\eta^2)=0$. When $c=0$, $R_\infty=-1/b$ and $R_\pm=b$:
this is the case that the selfdual Einstein metric is a bi-axial Bianchi IX
metric, i.e., of local cohomogeneity one under $\Un(2)$, and the distinguished
Einstein--Weyl quotient (by the centre of $\Un(2)$) is $S^3$~\cite{Ped:emm}.

In general, one checks that $R_\infty<R_\pm$ for all $\theta$.  Hence there
are three domains of definition.
\begin{itemize}
\item $R\in(-\infty,R_\infty)$: here the metric has negative scalar
curvature and yields a complete metric on the ball $B^4$ with a conformal
infinity at $R=R_\infty$.
\item $R\in(R_\infty,R_\pm)$: the metric still has negative scalar
curvature, but has an unremovable singularity at $R=R_\pm$.
\item $R\in(R_\pm,\infty)$: the metric now has positive scalar curvature,
again with an unremovable singularity at $R=R_\pm$.
\end{itemize}
The complete domain was found by Galicki~\cite{Gal:mcm}, although it is only
more recently that it has been noticed that these metrics are more general
than those of~\cite{Ped:emm}.

If $b=0$ then there are two nontrivial cases: for $c>1$, there is a conformal
infinity at $\sin\theta=-1/c$ and no singularity, while for $c<1$ there is a
singularity at $\sin\theta=-c$ and no conformal infinity. These turn out to be
Bianchi VIII metrics, i.e., have local cohomogeneity one under $\GL(2,\R)$,
and the distinguished Einstein--Weyl quotient is $\cH^3$~\cite{CaPe:sdc}.  For
$c=1$ the metric is the Bergmann metric on $\C\cH^2$.  We illustrate this
discussion with a diagram (Figure 2) of the $(b,c)$-plane next to which we
give a heuristic picture of the behaviour of $F$ on the branched double cover
of the hyperbolic plane by shading the domain over which the Einstein metric
has positive scalar curvature and indicating the conformal infinity in the
domain of negative scalar curvature.

\medbreak

\ifaddpics
\begin{center}
\includegraphics[width=.5\textwidth]{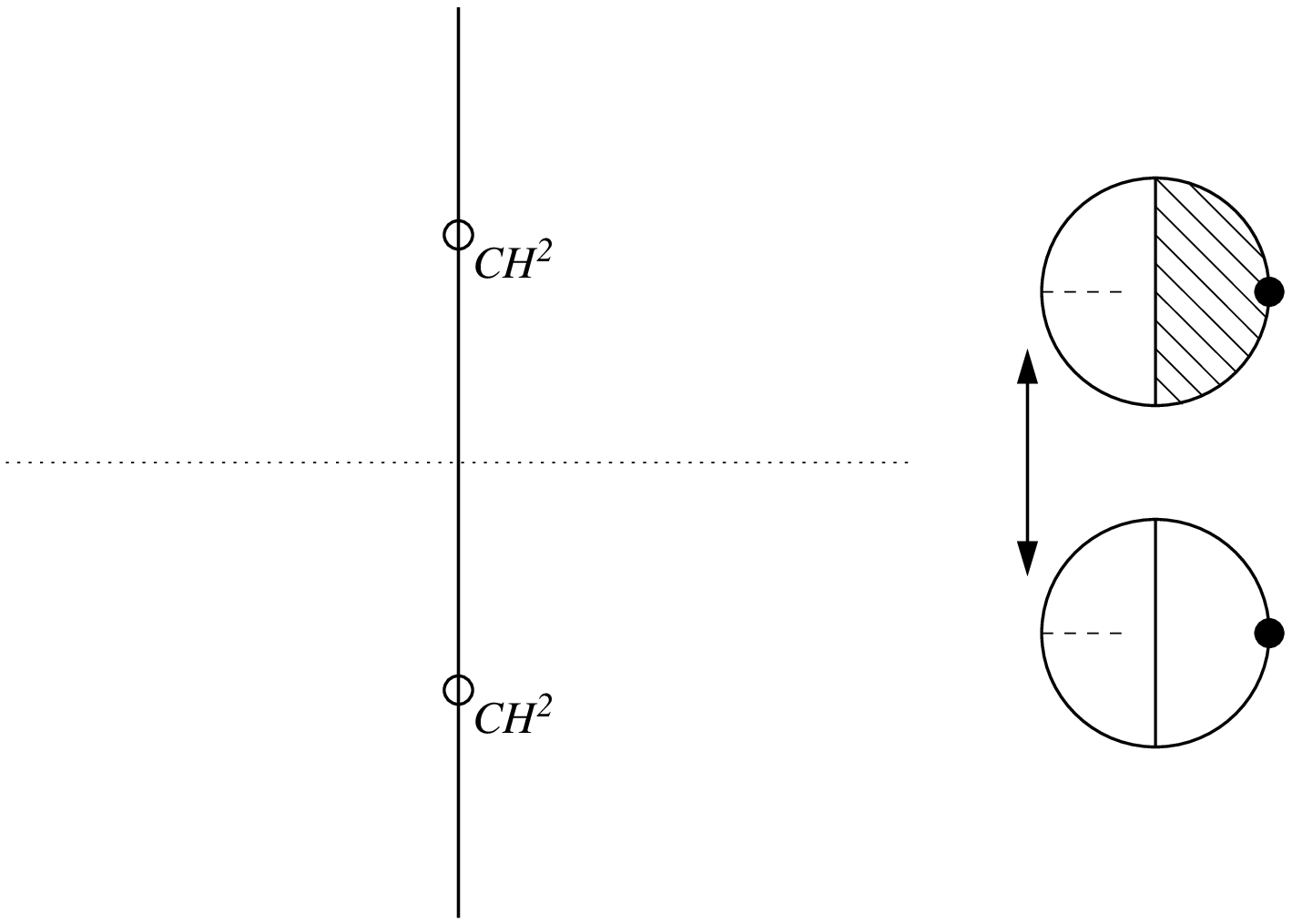}

\medbreak
Figure 2.
\end{center}
\fi

\bigbreak

For the Type II metrics the range of $R$ is $(1,\infty)$ and $F$ is globally
defined on the hyperbolic plane. However the moduli space has a much richer
structure and the analysis is more involved. The homothety freedom means that
the Einstein metrics are parameterized by a point $[a,b+c,b-c]$ in
$\RP2\punc{[1,0,0],[0,1,0],[0,0,1]}$. The lines joining the points $[1,0,0]$,
$[0,1,0]$ and $[0,0,1]$ represent dipole solutions, so the true moduli space
is obtained by removing these lines and taking the quotient by the permutation
group $\Sym_3$ of the coordinates. For convenience we shall only remove the
line $a=0$, so that we can set $a=1$ and use inhomogeneous coordinates $(b,c)$
as in the Type I case. On the lines $b=\pm c$ the selfdual Einstein metric is
the hyperbolic metric for $b<0$ and the spherical metric for $b>0$.

First note that $(b,c)=(1,0)$ (the fixed point $[1,1,1]$) gives the
Fubini--Study metric on $\CP2$, whereas the points $(0,1)$, $(0,-1)$ and
$(-1,0)$ yield the Bergmann metric on $\C\cH^2$ (cf.~\cite{GaLa:qro,Ped:emm}).
Along the lines joining these four points, we have bi-axial Bianchi metrics
with distinguished quotient $\cH^3$: along the lines joining $(1,0)$ to the
others, the metric is Bianchi IX, whereas on the lines between $(0,1)$,
$(0,-1)$ and $(-1,0)$, the metric is Bianchi VIII.

As in the Type I case, the zero-sets of $F$ and $\tfrac14 F^2-\rho^2(F_\rho^2
+ F_\eta^2)$ do not meet. This is a matter of checking that there are no
simultaneous solutions of
\begin{equation}\label{crux}\begin{split}
bR+cS+1&=0\\
b R + b^2 - cS - c^2&=0
\end{split}\end{equation}
with $S=\sin\theta\in(-1,1)$, $R>1$. The solution of~\eqref{crux} for
$b,c\neq0$ is $R=-(1 + b^2 - c^2)/2b, S=-(1 - b^2 + c^2)/2c$ which satisfies
$b^2(R^2-1)=c^2(S^2-1)$.  Hence we see that we cannot have $S\in[-1,1]$,
$R\in[1,\infty]$ unless $(b,c)$ lies on one of the Bianchi VIII lines $b=0$,
$b+c+1=0$ or $b-c+1=0$, in which case there are solutions $R=\infty$,
$(R,S)=(1,-1)$ and $(R,S)=(1,1)$ respectively, i.e., at the three marked
points at infinity.

Analysing the equations in~\eqref{crux} separately, we can determine for
which $(b,c)$ the selfdual Einstein metric has positive and/or negative
scalar curvature domains, and whether there is a conformal infinity.
The Bianchi VIII lines $b=0$, $b\pm c+1=0$ and the dipole lines $b=\pm c$
divide the $(b,c)$ plane into $1+3+3+6=13$ regions of four different types
leading to the following picture of the \mbox{(pre-)moduli} space (Figure 3).

\medbreak

\ifaddpics
\begin{center}
\includegraphics[width=.6\textwidth]{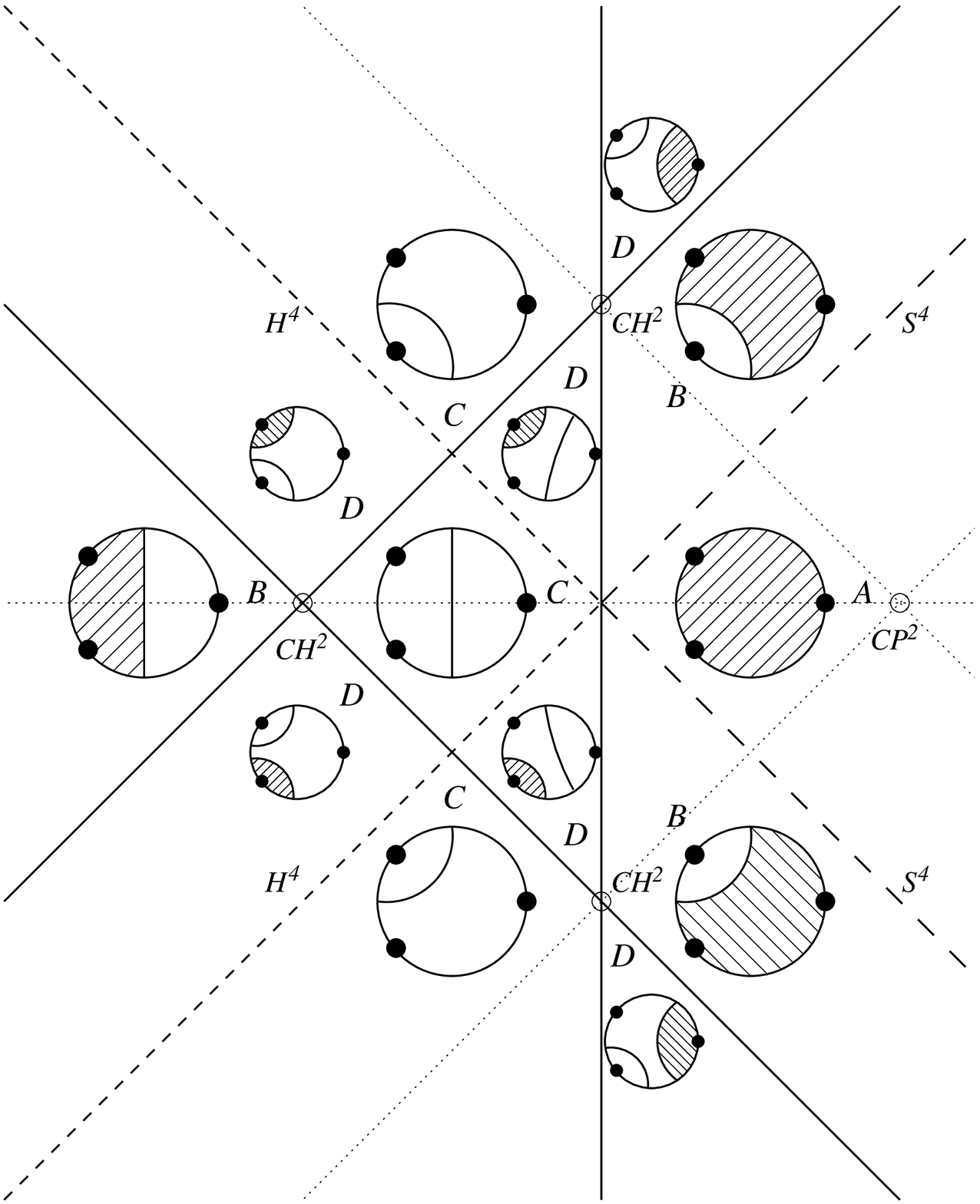}

\medbreak
Figure 3.
\end{center}
\fi

\bigbreak

In this figure the Bianchi VIII lines are solid, the Bianchi IX lines are
dotted (and do not bound regions), while the dipole lines are dashed. We label
the type of the region by $A,B,C,D$ and sketch the topology of the zero sets
in the hyperbolic disc as before. In these regions of the moduli space the
geometry of the selfdual Einstein metric is as follows.

\begin{enumerate}\renewcommand{\labelenumi}{\textup{(\Alph{enumi})}}
\item The metric has positive scalar curvature and no singularities for
$R\in(1,\infty)$. When $(b,c)=(1,0)$ it is the Fubini--Study metric, while for
other rational values the metric may be compactified on a weighted projective
space $\CP{[p,q,r]}$, and the selfdual Einstein metrics are
Hermitian~\cite{ApGa:seh,Bry:bkm,GaLa:qro}.
\item The metric has a positive and a negative scalar curvature domain
separated by an unremovable singularity.
\item The metric has two domains of negative scalar curvature separated by a
conformal infinity. On one side of the conformal infinity the metric is
complete on $B^4$. For rational parameter values, the metric on the other side
yields a complete metric on $\cO(n)\to\CP1$ or an orbifold generalization of
this~\cite{Gal:nmc,Hit:tem,LeBr:pac,Ped:emm}.
\item The metric has two domains of negative scalar curvature and one of
positive scalar curvature. On one side of the conformal infinity we obtain a
complete metric on $B^4$. This is similar to the behaviour in the Type I case.
\end{enumerate}

It is also fairly clear from Figure 3 what happens as $(b,c)$ passes from
one region to another. Along the dipole lines $b=\pm c$, a bubble of positive
or negative scalar curvature appears or disappears at one of the marked
points, whereas along the Bianchi VIII lines, either the singularity or the
conformal infinity passes through one of the marked points, creating or
destroying a domain of positive or negative scalar curvature as it does.

The $R,S$ coordinates give a reasonably simple formula for the metric.
\begin{align*}
g&=\frac{b^2-c^2+a (bR - c S) }{{{(a+b R+c S)}^2}}
\Bigl(\frac{dR^2}{R^2-1}+\frac{dS^2}{1-S^2}\Bigr)\\
&+\frac1{ (a+b R+c S)^2 (b^2-c^2 +a (b R - c S) ) (R^2-S^2)}\\
&\quad*\Bigl((R^2-1) (1-S^2)\bigl((b R-c S) d\phi+(cR-b S)d\psi\bigr)^2\\
&\qquad+\bigl((b (R^2-1) S+c (1- S^2)R) d\phi+(c (R^2-1) S+b(1-S^2)R
+a (R^2-S^2)) d\psi\bigr)^2\Bigr).
\end{align*}
In particular, the metric is rational. This is not surprising in view of the
work of Apostolov--Gauduchon~\cite{ApGa:seh} and Bryant~\cite{Bry:bkm}: the
$3$-pole metrics are all Hermitian, and are given explicitly in terms of a
fourth order polynomial $p(y)$ whose roots sum to zero. Despite the
superficial resemblence of this formula to~\cite[section 4.3.2]{Bry:bkm} (see
also~\cite{ACG:wsdk}), the precise relationship is rather complicated. By
computing the selfdual Weyl curvature of the Einstein metric, we find that the
selfdual K\"ahler metric is (a constant multiple of) $(a+bR+cS)^2g/(b^2-c^2+a
(bR - c S))^2$ with K\"ahler form
\begin{equation*}
\frac{1}{(b^2-c^2 +a (b R - c S))^2}\Bigl(
(d\phi\wedge (b\,dR-c\,dS)+d\psi\wedge \bigl((c+a S)dR-(b+a R)dS)\bigr)\Bigr).
\end{equation*}
The momentum maps of $\del_\phi$ and $\del_\psi$ are therefore (up to an
affine transformation)
\begin{equation*}
\frac{b R-c S}{b^2-c^2 +a (b R - c S)}\quad\text{and}\quad
\frac{c R - b S}{b^2-c^2 +a (b R - c S)}.
\end{equation*}
According to~\cite{Bry:bkm}, these must be affine linear combinations of the
trace $u_1=y_1+y_2$ and the pfaffian $u_2=y_1y_2$ of the normalized Ricci form
of the K\"ahler metric, and by~\cite{ApGa:seh} $u_1^2$ is the conformal factor
from the Einstein metric to the K\"ahler metric.  For generic $a,b,c$, this
allows us to put the K\"ahler metric into the form of~\cite{Bry:bkm,ACG:wsdk}:
\begin{equation*}
(y_1-y_2)\Bigl( \frac{dy_1^2}{p(y_1)}-\frac{dy_2^2}{p(y_2)}\Bigr)
+\frac1{y_1-y_2} \bigl ( p (y_1) (dt_1 + y_2 dt_2) ^2
- p (y_2) (dt_1 + y_1 dt_2) ^2 \bigr)
\end{equation*}
where $p(y)=(y-2 a b-{b^2}+{c^2}) (y+2 a b-{b^2}+{c^2}) (y-2 a c+{b^2}-{c^2})
(y+2 a c+{b^2}-{c^2})$. We interpret this formula abstractly by noting that
three twistors $\varphi_1,\varphi_2,\varphi_3$ give rise to three
$\SL(\W)$-invariants, namely the determinants $z_1=\eps(\varphi_2,\varphi_3)$,
$z_2=\eps(\varphi_3,\varphi_1)$ and $z_3=\eps(\varphi_1,\varphi_2)$. Up to an
overall sign, the roots of the polynomial $p(y)$ are then
\begin{align*}
r_0&=\tfrac12(z_1+z_2+z_3), &r_1 &= \tfrac12(z_1-z_2-z_3),\\
r_2&=\tfrac12(-z_1+z_2-z_3),&r_3 &= \tfrac12(-z_1-z_2+z_3).
\end{align*}
Hence we have related the generic Type II $3$-pole Einstein metrics to the
Case 4 K\"ahler metrics of~\cite{Bry:bkm}. In a similar way, the generic Type
I $3$-pole Einstein metrics are related to Bryant's Case 1 K\"ahler metrics.
The extra lines in our $3$-pole moduli space are cohomogeneity one metrics,
which are treated separately in~\cite{ApGa:seh,Bry:bkm} because the trace and
the pfaffian of the normalized Ricci form are not independent.  On the other
hand Bryant's Case 2 and Case 3 metrics are not covered by the $3$-pole
metrics, because $p(y)$ then has repeated roots. These correspond to
quaternion-\kahl/ quotients of $\cH\HQ^2$ and $\cH\HQ^{1,1}$ by a
non-semisimple $S^1$ action. We could obtain them from a limiting process
in which a dipole becomes infinitesimal.

Let us end by remarking that it is straightforward to obtain many complete
selfdual Einstein metrics with $T^2$ symmetry on $B^4$. Such metrics arise
when there is a domain of negative scalar curvature surrounding a single
marked point and bounded by a conformal infinity. Starting with a known
example we can deform $F$ slightly by adding additional monopole solutions at
points on the other side of the conformal infinity. The zero-set of $F$
deforms smoothly and so the metric stays complete until the conformal infinity
hits a fixed point. This argument yields not only the $3$-pole solutions, but
$m$-pole solutions for any $m>2$.  Hence the moduli space of smooth and
complete torus-symmetric selfdual Einstein metrics on the ball is infinite
dimensional, cf.~\cite{LeBr:cqk}.

We can also obtain infinite dimensional families of smooth complete metrics
with other topologies, but we postpone the discussion of these examples to
another occasion.

%
\newcommand{\bauth}[1]{\mbox{#1}} \newcommand{\bart}[1]{\textit{#1}}
\newcommand{\bjourn}[4]{#1\ifx{}{#2}\else{ \textbf{#2}}\fi{ (#4)}}
\newcommand{\bbook}[1]{\textsl{#1}}
\newcommand{\bseries}[2]{#1\ifx{}{#2}\else{ \textbf{#2}}\fi}
\newcommand{\bpp}[1]{#1} \newcommand{\bdate}[1]{ (#1)} \def\band/{and}
\newif\ifbibtex
\ifbibtex
\bibliographystyle{genbib}
\bibliography{papers}
\else

\fi

\end{document}